\documentclass[hidelinks,onefignum,onetabnum]{siamart220329}

\usepackage{lipsum}
\usepackage{amsfonts}
\usepackage{graphicx}
\usepackage{epstopdf}
\usepackage{algorithmic}
\usepackage{booktabs}
\usepackage{cleveref}
\usepackage{tikz}
\usepackage{bm}
\usepackage{caption}
\usepackage{float}
\usepackage{subcaption}

\ifpdf
  \DeclareGraphicsExtensions{.eps,.pdf,.png,.jpg}
\else
  \DeclareGraphicsExtensions{.eps}
\fi

\newsiamremark{remark}{Remark}
\newsiamremark{hypothesis}{Hypothesis}
\crefname{hypothesis}{Hypothesis}{Hypotheses}
\newsiamthm{claim}{Claim} 
\newsiamremark{example}{Example}

\newcommand{\xx}{\bm{x}}
\newcommand{\uu}{\mathbf{u}}
\newcommand{\nn}{\mathbf{n}}
\newcommand{\ff}{\mathbf{f}}

\newcommand{\JJ}{\mathbf{J}}

\newcommand{\R}{\mathbb{R}} 

\newcommand{\sgn}{\pm}
\newcommand{\dx}{\, \mathrm{d} \xx}
\newcommand{\ds}{\, \mathrm{d} s}

\headers{A scalable solver for ionic electrodiffusion}{P.~Benedusi et al}
\title{Scalable approximation and solvers for ionic electrodiffusion in cellular geometries}

\author{
  Pietro Benedusi\thanks{Department for Numerical Analysis and Scientific Computing, Simula Research Laboratory, Oslo, Norway (\email{benedp@simula.no}, \email{ada@simula.no}, \email{hherlyng@simula.no}, \email{meg@simula.no})} 
  \and Ada J.~Ellingsrud\footnotemark[1]
  \and Halvor Herlyng\footnotemark[1]
  \and Marie E.~Rognes\footnotemark[1] 
}

\usepackage{amsopn}

\ifpdf
\hypersetup{
  pdftitle={KNP-EMI Solvers},
  pdfauthor={P.~Benedusi et al}
}
\fi

\begin{document}

\maketitle

\begin{abstract}
 The activity and dynamics of excitable cells are fundamentally regulated and moderated by extracellular and intracellular ion concentrations and their electric potentials. The increasing availability of dense reconstructions of excitable tissue at extreme geometric detail pose a new and clear scientific computing challenge for computational modelling of ion dynamics and transport. In this paper, we design, develop and evaluate a scalable numerical algorithm for solving the time-dependent and nonlinear KNP-EMI equations describing ionic electrodiffusion for excitable cells with an explicit geometric representation of intracellular and extracellular compartments and interior interfaces. We also introduce and specify a set of model scenarios of increasing complexity suitable for benchmarking. Our solution strategy is based on an implicit-explicit discretization and linearization in time, a mixed finite element discretization of ion concentrations and electric potentials in intracellular and extracellular domains, and an algebraic multigrid-based, inexact block-diagonal preconditioner for GMRES. Numerical experiments with up to $10^8$ unknowns per time step and up to 256 cores demonstrate that this solution strategy is robust and scalable with respect to the problem size, time discretization and number of cores. 
\end{abstract}



\section{Introduction}

Our brains are composed of intertangled tissue consisting of neurons,
glial cells and interstitial space, penetrated by blood vessels. Ions,
such as potassium (K$^+$), sodium (Na$^+$), calcium (Ca$^{2+}$) and
chloride (Cl$^{-}$), move within and between the intracellular and
extracellular compartments in a highly regulated manner. These ions,
their movement and concentration differences, are fundamental for the
function and well-being of the brain~\cite{sterratt2023principles,
  rasmussen2020interstitial}. For instance, the brain's electrical
signals fundamentally rely on action potentials induced by rapid
neuronal influx of Na$^+$ and efflux of
K$^+$~\cite{sterratt2023principles}. Moreover, ion concentrations and
importantly ion concentration gradients may modulate brain signalling,
regulate brain volume, and control brain
states~\cite{nicholson1978calcium, aitken1986sources,
  rasmussen2020interstitial}. Notably, these physiological processes
are only partially understood and currently attract substantial
interest from the neurosciences~\cite{lagache2019electrodiffusion,
  rasmussen2020interstitial, yang2022carbon, armbruster2022neuronal,
  untiet2023astrocytic, nordentoft2023local,
  dietz2023local}. High-fidelity in-silico studies would offer an
innovative avenue of investigation with significant scientific
potential.

Electrical and chemical activity in excitable tissue are naturally
modelled via coupled partial differential equations (PDEs) describing
ionic electrodiffusion~\cite{mori2009numerical,
  niederer2013regulation, mori2015multidomain, ellingsrud2021accurate,
  ellingsrud2020finite}. These models go beyond the classical
electrophysiology models, such as the monodomain or bidomain
equations~\cite{sundnes2006computational} and their cellular-level
counterparts~\cite{agudelo2013computationally, tveito2017cell}, by
describing spatial and temporal variations of ion concentrations in
addition to the electric potentials. Within this context, we here consider
explicit geometric representations of the intracellular and
extracellular compartments and their joint interface representing the
cell membrane~\cite{mori2009numerical, yao2017numerical,
  ellingsrud2020finite}; i.e., the so-called KNP-EMI model
(Kirchhoff-Nernst-Planck, extracellular-membrane-intracellular).  In
their seminal research~\cite{mori2009numerical}, Mori and Peskin
introduced and analyzed a finite volume-based spatial discretization
of these cell-based ionic electrodiffusion equations, with numerical
experiments in idealized two- and three-dimensional
geometries. Similar equations and approximation challenges are also
encountered in connection with ionic transport in lithium-ion
batteries~\cite{kespe2019three}. In our previous
work~\cite{ellingsrud2020finite}, we introduced a mixed-dimensional
mortar finite element method for these equations, again with numerical
experiments limited to idealized geometries, and with little attention
to computational performance.

Dense reconstructions and segmentations of brain tissue are
increasingly becoming openly available~\cite{motta2019dense,
  zisis2021architecture, microns2021functional, lee20203d}. These
imaging-based representations reveal an extraordinary geometric
complexity that is very different from the more structured layers of
cardiac tissue~\cite{boron2016medical, jaeger2022arrhythmogenic}. This
complexity places clear demands on the scalability of numerical
solution algorithms. Solvers for the EMI (or cell-by-cell)
electrophysiology equations, which can be viewed as a subsystem of the
ionic electrodiffusion equations, have been studied quite intensively
in recent years~\cite{agudelo2013computationally, tveito2017cell,
  benedusi2024EMI}, including finite difference based
methods~\cite{jaeger2021efficient}, mixed finite element methods with
and without Lagrange multipliers~\cite{tveito2017cell}, boundary
element methods~\cite{rosilho2024boundary}, cut finite element
methods~\cite{berre2023cut}, and finite volume
methods~\cite{xylouris2010three}. Another key topic is scalable
preconditioning~\cite{huynh2023convergence,benedusi2024EMI}, and we
highlight the scalable solution strategy for electrodiffusion (without
interior interfaces) using a discontinuous Galerkin discretization
proposed by Roy et al~\cite{roy2023scalable}. Parallel scalability may
also be approached via domain decomposition techniques
\cite{huynh2023convergence, jaeger2021efficient}, though the
transition to complex geometries pose additional challenges. In fact,
most studies consider highly idealized
settings~\cite{mori2009numerical, ellingsrud2020finite}, structured
cell patterns as encountered in cardiac excitable
tissue~\cite{tveito2017cell, jaeger2022arrhythmogenic}, or simplified
neuronal
geometries~\cite{vossen2007modeling,xylouris2015three,shirinpour2021multi,buccino2021improving,berre2023cut, farina2021cut}.

In this paper, we address the challenge of how to design scalable
solution algorithms for the cell-based ionic electrodiffusion
(KNP-EMI) equations with high geometric complexity and physiological
membrane mechanisms. To this end, we present a lightweight
single-dimensional finite element formulation with a scalable and
efficient monolithic preconditioning strategy, vetted by a series of
numerical experiments with realistic model parameters and
geometries. The approximation scheme and solution strategy yield
accurate solutions in a bounded (low) number of Krylov
iterations. Near-optimal parallel scalability allows for rapidly
obtaining solutions to large KNP-EMI problems using high performance
computing systems, thus bridging the technology gap between dense
tissue reconstructions and ionic electrodiffusion simulations.

This manuscript is organized as follows. In Section~\ref{sec:model},
we present the KNP-EMI equations and describe active and passive
membrane dynamics involved in the interface coupling. In
Section~\ref{sec:disc}, we introduce a weak formulation of the
continuous problem together with discretizations in time and space. In
Section~\ref{sec:structure:solver}, we examine the block-structure and
properties of the resulting discrete systems, and present a tailored
solution strategy. In Section~\ref{sec::models} we present four model
scenarios in idealized and image-based geometries of increasing
complexity defining a benchmark suite for KNP-EMI
solvers. Section~\ref{sec::exp} reports on numerical experiments for
the four scenarios, investigating robustness, efficiency, and parallel
performance of the proposed strategy. Finally,
Section~\ref{sec:conclusions} provides some concluding remarks and
outlook.

\section{A model for ionic electrodiffusion in cellular spaces}
\label{sec:model}

Following~\cite{ellingsrud2020finite}, we consider the domain $\Omega = \Omega_i \cup \Omega_e \cup \Gamma \subset \R^d$ for $d \in \{2,3\}$ with $\Omega_i$ and $\Omega_e$ representing the intra- and extracellular domains (ICS, ECS), respectively, and $\Gamma = \partial\Omega_i$ the cellular membrane(s). In general, the intracellular domain can be composed of $N_{\rm cell}\in \mathbb{N}$ disjoint domains $\Omega_i = \bigcup_{j=1}^{N_{\rm cell}}\Omega_{i,j}$, with $\Omega_{i,j}$ representing the $j$th (biological) cell with membrane $\Gamma_j=\partial\Omega_{i,j}$. We consider a set $K$ of ionic species, for example $K = \{\text{Na}^+, \text{K}^+, \text{Cl}^-\}$,  with cardinality $|K|$. For each ionic species $k \in K$ and physical region $r\in\{i,e\}$, we model the evolution of ionic concentrations $[k]_r: \overline{\Omega}_r \times [0,T] \to \R$, as well as electric potentials $\phi_r:\overline{\Omega}_r \times[0,T] \to \R$.

\subsection{Conservation equations}

The assumption of conservation of ions for each region $r \in \{i, e\}$ yields the following time-dependent partial differential equation for each ionic species $k \in K$ and for $t\in (0,T)$ and $\xx\in\Omega_r$:
\begin{equation}
  \partial_t [k]_r(\xx,t) + \nabla \cdot \JJ_r^k(\xx,t) = f_r^k(\xx,t),
  \label{eq:conservation-ions}
\end{equation}
where $\JJ_r^k:\overline{\Omega}_r\times[0,T]\to\R^d$ is the ion flux density. We consider a Nernst--Planck relation describing diffusion and drift in the electric field induced by the potential $\phi_r$:
\begin{equation}
  \JJ_r^k = \JJ_{r,\mathrm{diff}}^k + \JJ_{r,\mathrm{drift}}^k  \\
   = - D_r^k \nabla [k]_r - \frac{D_r^k z_k}{\psi} [k]_r \nabla \phi_r,
  \label{eq:flux}
\end{equation}
where $D_r^k\in\R_+$ is the effective diffusion coefficient of ionic species $k$ in region $r$, $z_k$ is the valence of ion $k$, and $\psi = R T F^{-1}$ is composed of the gas constant $R$, the absolute temperature $T$, and the Faraday constant $F$. Finally, $f_r^k:\Omega_r\times[0,T]\to\R$ represent source terms. For consistency with the Kirchhoff-Nernst-Planck (KNP) assumption below, $f_r^k$ should maintain electroneutrality in the sense that
\begin{equation*}
  \sum_{k\in K} z_kf_r^k(\xx,t)=0  \quad \text{ for } r \in \{i,e\}.
\end{equation*}
Equation \eqref{eq:conservation-ions} can be interpreted as an advection-diffusion equation with advection velocity given by the electric field $\nabla \phi_r$. Since $\phi_r$ is itself an unknown, it is not trivial to classify \eqref{eq:conservation-ions} as a diffusion- or advection-dominated problem a-priori. We note that \eqref{eq:conservation-ions} ignores convective effects due to movement of the domain or medium; thus assuming a medium at rest (see e.g.~\cite{ellingsrud2021accurate, saetra2023neural} for contrast).

We introduce the KNP assumption of bulk electroneutrality for each region $r \in \{i, e\}$, for $t \in [0,T]$ and $\xx \in \Omega_r$:
\begin{equation*}
  \sum_{k\in K} z_k [k]_r (\xx,t) = 0,
\end{equation*}
which, differentiated with respect to time and combined with \eqref{eq:conservation-ions}, yields
\begin{equation}
\sum_{k\in K}z_k \nabla\cdot\JJ_r^k(\xx,t) = 0 .
  \label{eq:knp}
\end{equation}
The KNP assumption can be compared to the Poisson-Nernst-Planck (PNP) assumption. The latter would explicitly account for also the (rapid) charge relaxation dynamics, more relevant at the nanoscale~\cite{dickinson2011electroneutrality, solbraa2018kirchhoff}.  Note that the ion fluxes nonlinearly couple the ion concentrations $[k]_r$ via the potential $\phi_r$. Given the $|K| + 1$ unknown scalar fields in each of the two regions, the KNP-EMI model thus consists of two systems of $|K| + 1$ partial differential equations ($|K|$ parabolic, one elliptic). These two systems are coupled via EMI-type interface conditions.

\subsection{Interface conditions}

In this section, we introduce interface conditions describing currents and dynamics across the cellular membrane(s) $\Gamma$. We first define the membrane potential $\phi_M : \Gamma \times (0,T] \to \R$ as the jump in the electric potential:
\begin{equation}
  \phi_M(\xx,t) = \phi_i(\xx,t) - \phi_e(\xx,t)  \quad \text{for } \xx \in \Gamma,\, t > 0.
   \label{eq:membrane-potential}
\end{equation}
We next assume that the total ionic current density $I_M:\Gamma\times(0,T]\to\R$ is continuous across $\Gamma$ such that for $t > 0$ and $\xx \in \Gamma$,
\begin{equation}\label{eq::I_M1}
  I_M(\xx,t) = F \sum_{k \in K} z_k \JJ_i^k(\xx,t) \cdot \nn_i(\xx) = - F \sum_{k \in K} z_k \JJ_e^k(\xx,t) \cdot \nn_e(\xx),
\end{equation}
where $\nn_r$ denotes the normal on the interface pointing out from $\Omega_r$. We further assume that $I_M$ consists of two components: a total ionic channel current $I_{\rm ch}$ and a capacitive current $I_{\rm cap}$. Both of these currents have contributions from each ionic species: 
\begin{equation}\label{eq::I_M2}
I_M = I_{\rm ch} + I_{\rm cap} = \sum_{k \in K} I_{\rm ch}^k + \sum_{k \in K} I_{\rm cap}^k .
\end{equation}
The ionic channel currents $I_{\rm ch}^k$ will be subject to modelling, as described in Section~\ref{sec::ion}, while we consider the capacitor equation for the total capacitive current
\begin{equation}
I_{\rm cap} = C_m \partial_t \phi_M ,
\end{equation}
where $C_m\in\R_+$ is the membrane capacitance. 

As in~\cite[Section 2.1.3]{ellingsrud2020finite}, we derive an expression for $\JJ_r^k \cdot \nn_r$ that will be of useful in the weak formulation of the KNP-EMI equations. For $k \in K$ and $r \in \{i, e\}$, we introduce the relation
\begin{equation}
  \label{eq::I_cap}
I_{\rm cap}^k = \alpha_r^kI_{\rm cap} = \alpha_r^k C_m \partial_t \phi_M,
\end{equation}
and define the ratio
\begin{equation}
  \alpha_r^k = \frac{D_r^kz_k^2[k]_r}{\sum_{\ell \in K}D_r^\ell z_\ell^2[\ell]_r} \in [0,1],
  \quad\text{with} \quad \sum_{k \in K} \alpha_r^k=1.
\end{equation}
Combining~\eqref{eq::I_M1}, \eqref{eq::I_M2} and \eqref{eq::I_cap}, we can express the normal fluxes corresponding to a specific ion $k\in K$, as a function of the ionic currents and the membrane potential:
\begin{equation}\label{eq::J_n}
\JJ_r^k\cdot \nn_r= \sgn_r\frac{I_{\rm ch}^k + I_{\rm cap}^k}{Fz_k}=\sgn_r\frac{I_{\rm ch}^k + \alpha_r^kC_m \partial_t \phi_M}{Fz_k},
\end{equation}
for $r \in \{i, e\}$, $\sgn_i=1$, and $\sgn_e=-1$.

\subsection{Initial and boundary conditions}
\label{sec::bc}
To close the KNP-EMI system, we impose initial conditions for $k\in K$ and $r\in\{i,e\}$:
\begin{align}
&[k]_r(\xx,0)=[k]_r^0(\xx), &&\text{for } \xx\in\Omega_r,\\ 
&\phi_M(\xx,0)=\phi_M^0(\xx), &&\text{for } \xx\in\Gamma,
\end{align}
with $[k]_r^0$ and $\phi_M^0$ representing initial ion concentration and an initial membrane stimulus, respectively. We also assume that the initial ion concentrations are prescribed such that bulk electroneutrality is satisfied in each region:
\begin{equation}
  \sum_{k \in K}z_k[k]_r^0(\xx) = 0.
\end{equation}
We impose homogeneous Neumann boundary conditions on the exterior boundary $\partial \Omega$
\begin{equation}\label{eq::bc}
  \JJ_r^k(\xx,t) \cdot \nn_r(\xx) = 0 \quad \text{for } \xx \in \partial\Omega.
\end{equation}
The electric potentials $\phi_r$ in \eqref{eq:flux} are only determined up to a (common) additive constant. Therefore, an additional constraint must be considered to enforce uniqueness. Previous works have for example considered a zero-mean condition on $\phi_e$~\cite{ellingsrud2020finite}. We return to this point in \Cref{rem:nullspace}.

\subsection{Ionic channels}\label{sec::ion}
In terms of ionic species, we focus on sodium, potassium, and chloride: $K =
\{\mathrm{Na}^{+}, \mathrm{K}^{+}, \mathrm{Cl}^{-}\}$, and two
modelling scenarios (active and passive dynamics) for the ionic
channel currents $I_{\rm ch}^k$.

\subsubsection{Active dynamics: Hodgkin-Huxley model}
\label{sec:act_dyn}
To model the membrane dynamics of axons, we consider the Hodgkin-Huxley (HH) ionic model \cite{hodgkin1952quantitative} with the following expressions for the ionic currents: 
\begin{equation}
  \begin{aligned}
    I_{\rm ch}^{\rm Na} (\phi_M,[{\rm Na}^+]_i,[{\rm Na}^+]_e,\bm{w}) &= ({g}^{\rm Na}_{\rm stim} + {g}^{\rm Na}_{\rm leak}+\Bar{g}^{\rm Na}m^3h)(\phi_M - E_{\rm Na}),\\
    I_{\rm ch}^{\rm K} (\phi_M,[{\rm K}^+]_i,[{\rm K}^+]_e,\bm{w}) &= ({g}^{\rm K}_{\rm leak}+\Bar{g}^{\rm K}n^4)(\phi_M - E_{\rm K}),\\
    I_{\rm ch}^{\rm Cl} (\phi_M,[{\rm Cl}^-]_i,[{\rm Cl}^-]_e) &= {g}^{\rm Cl}_{\rm leak}(\phi_M - E_{\rm Cl}),
    \label{eq:HH_ionic_model}
  \end{aligned}
\end{equation}
given the ion-specific reversal potentials
\begin{equation}
  E_k=\frac{RT}{Fz_k}\ln{\frac{[k]_e}{[k]_i}}.
\end{equation}
In~\eqref{eq:HH_ionic_model}, $\Bar{g}^k,g_{\rm leak}^k\in\R_+$ are the maximum and leak conductivities, respectively. The Hodgkin-Huxley model also includes a set of time-dependent gating variables $\bm{w}(t) = (n(t), m(t), h(t)) \in [0,1]^3$ governed by the initial value problem
\begin{align}
  \label{eq::gating}
   \frac{\partial \bm{w}}{\partial t} &= \alpha_{\bm{w}}(\phi_M)(1-\bm{w})-\beta_{\bm{w}}(\phi_M)\bm{w}, \\
   \bm{w}(0) &= \bm{w}_0 = (n_0,m_0,h_0).
\end{align}
The coefficients $\alpha_w,\beta_w$ depend on the difference $\phi_M - \phi_{\rm rest}$, given a resting potential $\phi_{\rm rest}$, see~\cite{hodgkin1952quantitative}. The stimulus contribution ${g}^{\rm Na}_{\rm stim}: \Gamma \times [0,T] \to \R_+$ triggers the activation of the membrane dynamics.

  \subsubsection{Passive dynamics: Kir--Na/K}
  \label{sec:pass_dyn}
To model astrocyte and dendrite membrane dynamics, we consider the
following passive model (labeled Kir--Na/K) , which includes leak
channels, an inward-rectifying K channel, and a Na/K pump~\cite{halnes2013electrodiffusive}:
\begin{align*}
    I_{\rm ch}^{\rm Na}(\phi_M,[{\rm Na}^+]_i,[{\rm Na}^+]_e) &= ({g}^{\rm Na}_{\rm stim} + {g}^{\rm Na}_{\rm leak}) (\phi_M - E_\mathrm{Na}) + 3Fz_\mathrm{Na}\cdot j_\mathrm{pump}, \\
    I_{\rm ch}^{\rm K}(\phi_M,[{\rm K}^+]_i,[{\rm K}^+]_e) &= {g}^{\rm K}_{\rm leak}
     (\phi_M - E_\mathrm{K})f\mathrm{_{Kir}} - 2 F z_\mathrm{K} \cdot j_\mathrm{pump}, \\
    I_{\rm ch}^{\rm Cl}(\phi_M,[{\rm Cl}^-]_i,[{\rm Cl}^-]_e) &= {g}^{\rm Cl}_{\rm leak}(\phi_M - E_\mathrm{Cl}),
\end{align*}
where the Kir-function $f_{\mathrm{Kir}}$, controlling the inward-rectifying
K current, is given by:
\begin{equation*}
  f_{\mathrm{Kir}}([\mathrm{K}^+]_e, \Delta\phi_\mathrm{K}, \phi_M) =
   \frac{AB}{CD} \sqrt{\frac{\mathrm{[K^+]}_e}{\mathrm{[K^+]}_e^0} }, 
\end{equation*}
where
\begin{align*}
  A &= 1 + \exp(0.433), \qquad
  B = 1 + \exp(-(0.1186+E_\mathrm{K}^0)/0.0441), \\
  C &= 1 + \exp((\Delta\phi_\mathrm{K} + 0.0185)/0.0425),  \quad
  D = 1 + \exp(-(0.1186+\phi_\mathrm{M})/0.0441).
\end{align*}
Here, $\Delta\phi_\mathrm{K}=\phi_M - E_\mathrm{K}$, and
$E_\mathrm{K}^0$ is the K-reversal potential at $t=0$. Finally, the
pump flux density is given by:
\begin{equation*}
    j_\mathrm{pump} = \rho_\mathrm{pump} \left( \frac{\mathrm{[Na^+]}_i^{1.5}}{\mathrm{[Na^+]}_i^{1.5} + P_\mathrm{Na}^{1.5}} \right)
    \left( \frac{\mathrm{[K^+]}_e}{\mathrm{[K^+]}_e + P_\mathrm{K}} \right),
    \label{eq:pump}
\end{equation*}
where $\rho_\mathrm{pump},P_k\in\R_+$ are, respectively, the maximum
pump rate and a threshold for ion $k$. All physical parameters are
given in~\Cref{tab:physconst} (\Cref{sec::models}).

\section{Variational formulation and discretization}
\label{sec:disc}

\subsection{Single-dimensional formulation of the KNP-EMI problem}

We consider a so-called \textit{single-dimensional formulation} of the KNP-EMI problem, following the naming convention introduced in \cite{kuchta2021solving, tveito2021modeling}. In this formulation, the membrane current variable $I_M$ is eliminated and the problem is no longer mixed-dimensional, since all the unknowns are defined over $d$-dimensional domains. For the EMI problem, using the single-dimensional formulation results in a more compact algebraic structure with favorable properties (symmetry and positive definiteness)~\cite{benedusi2024EMI}. In the case of KNP-EMI, only positive definiteness is maintained.

Let $V_r^k = V_r^k(\Omega)$ and $W_r = W_r(\Omega)$ be Hilbert spaces of sufficiently regular functions. For $r\in\{i,e\}$, assuming the solutions $[k]_r$ and $\phi_r$ to be sufficiently regular over $\Omega_r$, we multiply the PDE in \eqref{eq:conservation-ions} by test functions $v^k_r\in V_r^k$ and integrate over $\Omega_r$. After integration by parts and using the normal flux definition from \eqref{eq::J_n}, we obtain for $k\in K$ and $r \in \{r, e\}$:
\begin{equation}
  \label{eq::weak2}    
  \int_{\Omega_r} \partial_t[k]_r v_r^k - \JJ^k_r \cdot \nabla v_r^k  \dx
  \sgn_r \frac{1}{Fz_k}\int_{\Gamma}\left(I^k_{\text{ch}} + \alpha_r^kC_m\partial_t\phi_M\right)v_r^k \ds = \int_{\Omega_r}f_r^kv_r^k \dx
\end{equation}
Similarly, multiplying \eqref{eq:knp} by test functions $w_r \in W_r$ and integrating by parts, we obtain 
\begin{equation}
  \label{eq::weak3}
  - \sum_{k \in K} z_k\int_{\Omega_r}  \JJ^k_r\cdot\nabla w_r \dx
  \sgn_r \frac{1}{F}\int_{\Gamma}\left(I_{\text{ch}} + C_m\partial_t\phi_M\right)w_r \ds = 0.
\end{equation} 
We recall that $I_{\text{ch}}^k$ linearly depends on $\phi_M = \phi_i-\phi_e$ (cf. Section~\ref{sec::ion}) and non-linearly on $[k]_r$. We thus define the weak solution to the KNP-EMI problem as $\{[k]_r,\phi_r\}$ for $r\in\{i,e\}$ and $k \in K$ satisfying \eqref{eq::weak2}--\eqref{eq::weak3} for all test functions $v_r^k\in V_r^k$ and $w_r\in W_r$. 

\subsection{Time discretization}
We consider an implicit-explicit time discretization scheme, leading to a convenient linearization. 
We partition the time axes uniformly into $N_t$ intervals:
$$t_n = n\Delta t, \quad \text{with} \quad n = 0,\ldots,N_t, \quad \Delta t = T/N_t.$$
Given the initial conditions $[k]_r^0,\phi_M^0,\bm{w}_0$ defined in Sections~\ref{sec::bc}--\ref{sec::ion}, we consider the following first-order, implicit approximations at times $t_n$ for $n=1,\ldots,N_t$:
\begin{equation}
  \partial_t[k]_r \approx ([k]_r^n - [k]_r^{n-1})/\Delta t, \qquad  \partial_t\phi_M \approx (\phi_M^n - \phi_M^{n-1})/\Delta t,    
\end{equation}
with $[k]_r^n = [k]_r(\xx,t_n), \phi_r^n=\phi_r(\xx,t_n)$, and $f_r^{k,n}=f^k_r(\xx,t_n)$. The ion flux densities are approximated with an implicit-explicit scheme, while the coefficients $\alpha_r^{k}$ and ionic currents $I_{\rm ch}^k$ are treated explicitly:
$$ \JJ_r^{k}(\xx,t_n) \approx \JJ_r^{k,n} = - D_r^k \nabla [k]_r^n - \frac{D_r^k z_k}{\psi} [k]_r^{n-1} \nabla \phi_r^n,$$
    $$\alpha_r^{k}(\xx,t_n) \approx\alpha_r^{k,n-1}=\frac{D_r^kz_k^2[k]_r^{n-1}}{\sum_{\ell \in K}D_r^\ell z_\ell^2[\ell]_r^{n-1}}.$$
If gating variables $\bm{w}$ are present, we first solve the ODEs \eqref{eq::gating} numerically with $N_{\text{ode}}$ steps of the Rush-Larsen method with step size $\Delta t/N_{\text{ode}}$ to obtain $\bm{w}_{n}\approx\bm{w}(t_n)$, setting $\phi_M=\phi_M^{n-1}$ in \eqref{eq::gating}, given the initial condition $\bm{w}_{n-1}$. We then approximate the ionic currents by
\begin{equation}
  I^{k}_{\text{ch}}(\xx,t_n)\approx I^{k,n-1}_{\text{ch}}= I_{\rm ch}^k (\phi_M^{n-1},[k]_i^{n-1},[k]_e^{n-1},\bm{w}_n).
\end{equation}
Inserting the previous expressions into \eqref{eq::weak2}, we obtain the following linearized weak formulation:
\begin{multline}
  \label{eq::weak_disc} 
    \int_{\Omega_r} [k]_r^n v_r^k - \Delta t \,\JJ_r^{k,n} \cdot \nabla v_r^k \dx
    \sgn_rC_r^{k,n-1}\int_{\Gamma} \phi_M^nv_r^k \ds \\
    = \int_{\Omega_r} ([k]_r^{n-1} + f_r^{k,n}) v_r^k \dx  - \int_{\Gamma} g_r^{k,n-1}v_r^k \ds,
\end{multline}
for $k\in K$ and $n=1,\ldots,N_t$, having introduced the short-hand
\begin{equation}
  C_r^{k,n} = \frac{\alpha_r^{k,n}C_m}{Fz_k} \quad \text{and} \quad g_r^{k,n} = \sgn_r\frac{1}{Fz_k}(\Delta t I^{k,n}_{\text{ch}}-\alpha_r^{k,n}C_m\phi_M^n).
\end{equation}
Similarly, for \eqref{eq::weak3}, we obtain
\begin{equation}
  \label{eq::weak_disc2}
  - \Delta t \sum_{k \in K} z_k \int_{\Omega_r} \JJ^{k,n}_r\cdot\nabla w_r \dx
  \sgn_r \frac{C_m}{F} \int_{\Gamma}\phi_M^nw_r \ds = - \int_{\Gamma} h_r^{n-1} w_r \ds,
\end{equation}
with $h^n_r = \sgn_r(\Delta t I_{\text{ch}}^n + C_m\phi_M^n)/F$. 

In summary, the semi-discrete KNP-EMI problem reads as: for
$n=1, \ldots, N_t$, for $r \in \{i,e\}$ and $k \in K$, find $[k]_r^n,
\phi_r^n$ such that~\eqref{eq::weak_disc} and~\eqref{eq::weak_disc2}
hold for all $v_r^k \in V_r^k$ and $w_r \in W_r$.

\subsection{Spatial discretization}
We discretize $\Omega_i$ and $\Omega_e$ by conforming simplicial meshes $\mathcal{T}_i$ and $\mathcal{T}_e$, respectively, such that $\mathcal{T} = \mathcal{T}_i \cup \mathcal{T}_e$ forms a conforming simplicial tessellation of $\Omega$; in particular, that $\mathcal{T}_i$ and $\mathcal{T}_e$ match at $\Gamma$. For $\mathcal{T}_r$ and $r\in\{i,e\}$, we approximate the Hilbert spaces $V_r^k$ and $W_r$ by finite element spaces $V^k_{r,h}$ and $W_{r,h}$ of continuous piecewise polynomials of degree $p_{r,k}$ and $p_r$, for $[k]_r$ and $\phi_r$ respectively. The (fully) discrete KNP-EMI problem to be solved at each time step $n$ thus reads: find $[k]_{r,h}^n \in V^k_{r,h}$ and $\phi_{r,h}^n\in W_{r,h}$ such that \eqref{eq::weak_disc} and \eqref{eq::weak_disc2} hold for all $v_r^k\in V^k_{r,h}$ and for all $w_r\in W_{r,h}$.

In practice, we use the same polynomial order $p=p_r=p_{r,k}$ for all regions $r$ and unknowns, therefore using a single finite element space $V_{r,h}$ of dimension (number of degrees of freedom) $N_r$ with $V_{r,h}=V^k_{r,h}=W_{r,h}$ for all $k \in K$. We remark that extending the content of the current and subsequent sections to the general case, with multiple finite element spaces, is straightforward, but cumbersome in terms of notation. We use Lagrangian nodal basis functions $\{\varphi_{r,j}\}_{j=1}^{N_r}$ of order $p$ (omitting the order in the notation) so that 
$$V_{r,h}=\mathrm{span}\left(\{\varphi_{r,j}\}_{j=1}^{N_r}\right),$$
and the real coefficients $\{[k]_{r,j}^n\}_{j=1}^{N_r}$ and $\{\phi_{r,j}^n\}_{j=1}^{N_r}$ define the discrete approximations
\begin{equation}
  [k]_r^n(\xx)\approx [k]_{r,h}^n(\xx)=\sum_j^{N_r}[k]_{r,j}^n\varphi_{r,j}(\xx)
  \quad \text{and} \quad \phi_r^n(\xx)\approx \phi_{r,h}^n(\xx)=\sum_j^{N_r}\phi_{r,j}^n\varphi_{r,j}(\xx),
\end{equation}
for $\xx\in\Omega_r$ and $n=0,\ldots,N_t$.

The potential jump $\phi_M$ and gating variables $\bm{w}$ are approximated with the same Lagrangian elements restricted to the membrane. To illustrate, we have 
$$\phi_{M,h}^n(\xx)=\phi_{i,h}^n(\xx) - \phi_{e,h}^n(\xx)=\sum_{j \in J_{\Gamma}}\phi_{M,j}^n \varphi_{r,j}(\xx) \quad \text{for} \quad \xx\in\Gamma,$$ 
for a fixed but arbitrary $r$, and $J_{\Gamma} = \{1\leq j\leq
N_r:\text{supp}(\varphi_{r,j})\cap \Gamma \neq \emptyset\}$ being the
set of indices corresponding to $\Gamma$. When applicable (if gating
variables are present), we solve \eqref{eq::gating} point-wise in
space at the coordinates associated with the Lagrange degrees of
freedom and using the coefficient values $\phi_{M,j}$.

\section{Discrete algebraic structure and solver strategy}
\label{sec:structure:solver}

We now turn to examine the structure and properties of linear
operators associated with the discrete KNP-EMI equations, before
presenting a tailored iterative and preconditioned linear solution
strategy.

\subsection{Algebraic structure and matrix assembly}
\label{sec:structure}

We begin by defining mass $M$ and (weighted) stiffness matrices $A$
for the concentrations and potentials in each of the bulk regions
$\Omega_r$:
\begin{align*}    
    M_r & =\left[\int_{\Omega_r}\varphi_{r,j}(\xx)\varphi_{r,l}(\xx)\,\mathrm{d}\xx\right]_{j,l=1}^{N_r}\in\mathbb{R}^{N_r\times N_r}, \\
    A_r& =\left[\int_{\Omega_r}\nabla\varphi_{r,j}(\xx)\cdot\nabla\varphi_{r,l}(\xx)\,\mathrm{d}\xx\right]_{j,l=1}^{N_r}\in\mathbb{R}^{N_r\times N_r},\\
    A_r^{k,n} & =\left[\int_{\Omega_r}[k]_r^n(\xx)\nabla\varphi_{r,j}(\xx)\cdot\nabla\varphi_{r,l}(\xx)\,\mathrm{d}\xx\right]_{j,l=1}^{N_r}\in\mathbb{R}^{N_r\times N_r}, 
\end{align*}
We next define the membrane (mass) matrices:
\begin{align*}    
    M_{r,\Gamma} & =\left[\int_{\Gamma}\varphi_{r,j}(\xx)\varphi_{r,l}(\xx)\,\mathrm{d}s\right]_{j,l=1}^{N_r}\in\mathbb{R}^{N_r\times N_r}, \\
    M_{rq,\Gamma} & =\left[\int_{\Gamma}\varphi_{r,j}(\xx)\varphi_{q,l}(\xx)\,\mathrm{d}s\right]_{j,l=1}^{(N_r,N_q)}\in\mathbb{R}^{N_r\times N_q},
\end{align*}
with $q=\{i,e\}/r$.
Note that the membrane operators, denoted by the subscript $\Gamma$, have non-zero elements only for $(j,l)$ entries corresponding to basis functions with $\Gamma$ in their support. For example,
\begin{equation*}
  \Gamma\cap\text{supp}(\varphi_{r,j})\cap\text{supp}(\varphi_{r,l})\neq\emptyset \quad \Rightarrow \quad [M_{r,\Gamma}]_{j,l} \neq 0.
\end{equation*}
Therefore, membrane operators can be seen as low-rank perturbations of bulk ones. 

To recast the discrete KNP-EMI equations \eqref{eq::weak_disc}--\eqref{eq::weak_disc2} as a linear system, we define the non-symmetric matrix\footnote{With a minor abuse of notation, we use $k$ both to identify an ionic species, i.e. $k\in K$ and as an index for elements $K$, i.e. $k=1,...,|K|$.}
corresponding to $\Omega_r$, given $\tau_r^k=\Delta t D_r^k$ and $\tilde{\tau}_r^k=\tau_r^kz_k\psi^{-1}$:
\begin{equation*}
\mathcal{A}_r^n  =
\left[
\begin{array}{cccc|c}
M_r + \tau_r^1A_r & 0  & \cdots & 0  &  C_r^{1,n}M_{r,\Gamma} + \tilde{\tau}_r^1{A}_r^{1,n} \\
 0 & M_r +\tau_r^2A_r&  0 & \vdots &   C_r^{2,n}M_{r,\Gamma} + \tilde{\tau}_r^2A_r^{2,n}  \\
 \vdots & 0 & \ddots & 0 & \vdots  \\
 0 & \cdots & 0 & M_r +\tau_r^{|K|}A_r & C_r^{|K|,n}M_{r,\Gamma} + \tilde{\tau}_r^{|K|}A_r^{|K|,n} \\
  \midrule
  z_1\tau_r^1A_r & z_2\tau_r^2A_r & \cdots & z_{|K|}\tau_r^{|K|}A_r & C_mF^{-1}M_{r,\Gamma} + \Sigma_{k}z_k\tilde{\tau}_r^k A_r^{k,n}
\end{array}
\right],
\end{equation*}
with size $N_r(|K| + 1)$. More compactly,
\begin{equation*}
\mathcal{A}_r^n = 
\begin{bmatrix}    
\mathcal{A}_r^{11} & & & & \mathcal{A}_r^{1\phi} \\
& \mathcal{A}_r^{22} & & & \mathcal{A}_r^{2\phi} \\
& & \ddots & & \vdots \\
& & & \mathcal{A}_r^{|K||K|} & \mathcal{A}_r^{|K|\phi} \\
\mathcal{A}_r^{\phi 1} & \mathcal{A}_r^{\phi 2} & \cdots & \mathcal{A}_r^{\phi |K|} & \mathcal{A}_r^{\phi\phi}
 \end{bmatrix}_n
=
\begin{bmatrix}    
\mathcal{A}_r^{cc} & \mathcal{A}_r^{c\phi} \\
\mathcal{A}_r^{\phi c} & \mathcal{A}_r^{\phi\phi}
 \end{bmatrix}_n, 
\end{equation*}
with $\mathcal{A}_r^{cc}\in\R^{N_r|K|\times N_r|K|}$ and $\mathcal{A}_r^{\phi\phi}\in\R^{N_r\times N_r}$ the concentrations and potentials blocks, respectively. The last block-column of $\mathcal{A}_r^n$ is time-dependent and must be updated in each time step $t_n$.

The block diagonal structure of $\mathcal{A}_r^{cc}$ highlights that the different concentrations are coupled only via the potential. We define the following coupling matrices:
\begin{equation*}
\mathcal{B}_{rq}^n =
\left[
    \begin{array}{ccc|c}
0 & \cdots & 0  &  - C_r^{1,n}M_{rq,\Gamma} \\
0 & \cdots & 0  &  - C_r^{2,n}M_{rq,\Gamma}  \\
\vdots  & & \vdots  &  \vdots \\
0 & \cdots & 0  &  - C_r^{|K|,n}M_{rq,\Gamma}  \\
\midrule
0 & \cdots & 0  &  - C_mF^{-1}M_{rq,\Gamma}  \\
\end{array}\right], 
\end{equation*}
with size $N_r(|K| + 1)\times N_q(|K| + 1)$, or, more compactly, 
\begin{equation*}
\mathcal{B}_{rq}^n =
 \begin{bmatrix}    
0 & \mathcal{B}_{rq}^{c\phi} \\
0 & \mathcal{B}_{rq}^{\phi\phi}
 \end{bmatrix}_n,
 \end{equation*}
 with $\mathcal{B}_{rq}^{c\phi}\in\R^{N_r|K|\times N_q}, \mathcal{B}_{rq}^{\phi\phi}\in\R^{N_r\times N_q}$. The global matrix at time $t_n$, coupling variables in $\Omega_i$ and $\Omega_e$, thus reads
 \begin{equation*}    
 \mathcal{A}_n=
  \begin{bmatrix}    
    \mathcal{A}_i^n & \mathcal{B}_{ie}^n\\
    \mathcal{B}_{ei}^n & \mathcal{A}_e^n\\
  \end{bmatrix}
  =
  \begin{bmatrix}    
    \mathcal{A}_i^{cc} & \mathcal{A}_i^{c\phi} & 0 & \mathcal{B}_{ie}^{c\phi}\\
    \mathcal{A}_i^{\phi c} & \mathcal{A}_i^{\phi\phi} & 0 & \mathcal{B}_{ie}^{\phi\phi}\\
    0 & \mathcal{B}_{ei}^{c\phi} & \mathcal{A}_e^{cc} & \mathcal{A}_e^{c\phi}\\
    0 & \mathcal{B}_{ei}^{\phi\phi} & \mathcal{A}_e^{\phi c} & \mathcal{A}_e^{\phi\phi}
  \end{bmatrix}_n.
\end{equation*}
Let us remark that for multiple cells, i.e. $N_{\rm cell}>1$, $\mathcal{A}_i^n$ can be further decomposed in $N_{\rm cell}$ diagonal blocks coupled with ${A}_e^n$ via the coupling matrices $\mathcal{B}_{rq}^n$.

Finally, we define the solution vectors:
\begin{equation}
\bm{c}^{k,n}_r=[[k]_{r,1}^n,\ldots,[k]_{r,N_r}^n]\in\R^{N_r}, \qquad
\bm{c}_r^n=[\bm{c}^{1,n}_r,\bm{c}^{2,n}_r,\ldots,\bm{c}^{|K|,n}_r]\in\R^{N_r|K|},
\end{equation}
\begin{equation}
\bm{\phi}_r^n=[\phi_{r,1}^n,...,\phi_{r,N_r}^n]\in\R^{N_r},
\end{equation}
and recast the discrete variational problem \eqref{eq::weak_disc}--\eqref{eq::weak_disc2} as the following linear system of size $N=(N_i+N_e)(|K|+1)$, for $n=1,\ldots,N_t$: 
\begin{equation}\label{eq::lin_sys}
  \begin{bmatrix}    
    \mathcal{A}_i^{cc} & \mathcal{A}_i^{c\phi} & 0 & \mathcal{B}_{ie}^{c\phi}\\
    \mathcal{A}_i^{\phi c} & \mathcal{A}_i^{\phi\phi} & 0 & \mathcal{B}_{ie}^{\phi\phi}\\
    0 & \mathcal{B}_{ei}^{c\phi} & \mathcal{A}_e^{cc} & \mathcal{A}_e^{c\phi}\\
    0 & \mathcal{B}_{ei}^{\phi\phi} & \mathcal{A}_e^{\phi c} & \mathcal{A}_e^{\phi\phi} \\
  \end{bmatrix}_{n-1}
 \begin{bmatrix}
\bm{c}_i^n \\
\bm{\phi}_i^n \\
\bm{c}_e^n \\
\bm{\phi}_e^n 
\end{bmatrix}
=
\begin{bmatrix}
\bm{f}_i^{c,n}    \\
\bm{f}_i^{{\phi},n} \\
\bm{f}_e^{c,n}    \\
\bm{f}_e^{{\phi},n} 
\end{bmatrix},
\quad \Longleftrightarrow \quad \mathcal{A}_{n-1} \uu_n = \ff_{n},
\end{equation}
with $\bm{f}_r^{c,n} = \left[\bm{f}_r^{1,n},\bm{f}_r^{2,n},...,\bm{f}_r^{|K|,n}\right]^T\in\R^{N_r|K|}$ and
\begin{align*}
\bm{f}_r^{k,n} & = \left[\int_{\Omega_r}([k]_r^{n-1}(\xx)+f_r^{k,n}(\xx))\varphi_{r,j}(\xx)\,\mathrm{d}\xx -\int_{\Gamma}g_r^{k,n-1}(\xx)\varphi_{r,j}(\xx)\,\mathrm{d}s \right]_{j=1}^{N_r}\in\R^{N_r}, \\
\bm{f}_r^{\phi,n} & = \left[-\int_{\Gamma}h^{n-1}(\xx)\varphi_{r,j}(\xx)\,\mathrm{d}s\right]_{j=1}^{N_r}\in\R^{N_r}.
\end{align*}
Let us remark that the choice of $\mathcal{A}_{n-1}$, instead of $\mathcal{A}_{n}$, makes problem \eqref{eq::lin_sys} linear.

\begin{remark}
  The linear operator $\mathcal{A}_{n}$ in~\eqref{eq::lin_sys} admits
  a one-dimensional nullspace corresponding to the constant
 potential (see Section \ref{sec::bc}). There are
  several approaches to eliminate this nullspace, such as enforcing a zero-mean integral
  constraint via a Lagrange multiplier \cite{ellingsrud2020finite}. However, this approach
  leads to a dense row and column in $\mathcal{A}_{n}$,
  severely affecting parallel matrix assembly and iterative solver
  performance. Instead, we here, for the sake of simplicity and
  efficiency, introduce a point-wise Dirichlet condition on $\phi_e$:
  given an arbitrary $\xx_e \in \Omega_e$, we impose
  $\phi_e(\xx_e, t)=0$ for all $t\in(0,T]$. A third approach would be to prescribe the discrete nullspace in the iterative solution
    strategy. 
    \label{rem:nullspace}
\end{remark}

\subsection{Solution and preconditioning strategy}
\label{sec:solver}
To solve the non-symmetric linear system \eqref{eq::lin_sys}, we employ a preconditioned GMRES method using the block diagonal preconditioner 
\begin{equation}\label{eq::P_prec}
P_n  =
    \begin{bmatrix}    
\mathcal{A}_i^{cc} & 0 & 0 & 0\\
0 & \mathcal{A}_i^{\phi\phi} & 0 & 0\\
0 & 0 & \mathcal{A}_e^{cc} & 0\\
0 & 0 & 0 & \mathcal{A}_e^{\phi\phi}\\
\end{bmatrix}_n
=
\begin{bmatrix}    
\mathcal{A}_i^{11} &  &  & \\
&  \ddots & \\
& & \mathcal{A}_i^{\phi\phi} &  \\
& & & \mathcal{A}_e^{11} &  &  & \\
& & & &  \ddots & \\
& & & & & \mathcal{A}_e^{\phi\phi}\\
\end{bmatrix}_n\in\mathbb{R}^{N\times N},
\end{equation} 
composed of $2(|K| + 1)$ blocks, decoupling all concentrations and
potentials. In case of several disjoint intracellular
domains, the intracellular blocks can be further subdivided into block-diagonal matrices, with $N_{\rm cell}$ blocks each. By definition, $P_n$ is symmetric and
positive definite and corresponds to a simplified problem where
$\Omega_i$ and $\Omega_e$ are decoupled and the continuity equation
contains only diffusive contributions and no electric drifts, with
each concentration evolving independently. In relation to Remark~\ref{rem:nullspace}, the constant potential null space has to be eliminated also for $P_n$. 

In the simpler EMI setting, the elimination of off-diagonal blocks can be motivated by a spectral argument~\cite{benedusi2024EMI}, demonstrating that the membrane terms are \textit{zero-distributed}; i.e., represent a negligible contribution to the eigenvalues distribution for a mesh resolution increasing uniformly in space. As further motivation for this choice, numerical experiments suggest that diffusion dominates the drift induced by gradients in the electric potentials, at least for $k \in \{\mathrm{Na^+,K^+}\}$ at the scale and tissue types under consideration \cite{solbraa2018kirchhoff}. While the extension of the spectral theory to the KNP-EMI case has still to be developed, the numerical experiments presented in the subsequent sections show that $P_n$ is an effective and robust preconditioner. Let us note that, for physiologically relevant parameters, the condition number of $\mathcal{A}_n$ can exceed $10^{10}$, resulting in stagnation of GMRES if no preconditioning is adopted. 

In practice, the action of $P_n^{-1}$ (which can be efficiently computed block-wise with a preconditioned CG) can be approximated by a single algebraic multigrid (AMG) iteration applied monolithically~\cite{benedusi2021fast}. Moreover, we may use $P_0$ as a preconditioner for all time steps $n=1,\ldots, N_t$. This reduces assembly times, with little or no impact on GMRES convergence as compared to using the time-dependent preconditioner $P_n$. 

\section{Idealized and image-based brain tissue benchmark scenarios}
\label{sec::models}

\begin{table}
  \centering
  \small
  \begin{tabular}{clllll}
    \toprule
    & Geometry & Domains & Dynamics & Stimulus & $N$ (max) \\
    \midrule
    A & Idealized & ICS, ECS & HH & $g_{\rm stim}^{\rm Na} \neq 0$ & $8.7 \cdot 10^6$ \\
    B & Imaging & Astrocyte, ECS & Kir--Na/K& $f_e^{\rm Na},f_e^{\rm K} \neq 0$ & $5.3 \cdot 10^6$ \\
    C & Imaging & Dendrite, astrocytes, ECS & Kir--Na/K & $g_{\rm stim}^{\rm Na} \neq 0$ & $2.0 \cdot 10^6$ \\
    D & Imaging & Neurons, astrocyte, ECS & HH, Kir--Na/K & $g_{\rm stim}^{\rm Na} \neq 0$ & $1.0 \cdot 10^8$ \\
    \bottomrule
    \vspace{1em}
  \end{tabular}
  \caption{Overview of model scenarios, labeled Models A--D. The
    idealized Model A includes both two-dimensional and
    three-dimensional ($d = 2, 3$) setups. Models B--D employ
    different image-based 3D geometries, including one or more
    astrocytes or partial astrocytic processes, neurons and neuronal
    compartments (dendrites, axons, spines) at different scales in
    addition to the contiguous extracellular space (ECS). The ICS in
    Model A and the neurons in Model D are modelled using
    Hodgkin-Huxley (HH) membrane dynamics
    (Section~\ref{sec:act_dyn}). The astrocyte membranes in Models
    B--D and the dendrite membrane in Model C are modelled using the
    passive Kir--Na/K model (Section~\ref{sec:pass_dyn}). Each
    scenario including the stimulus are described in further detail in
    the respective model sections.}
  \label{tab:scenarios} 
\end{table}
In this section, we present four model scenarios for ionic
electrodiffusion in brain tissue, including neuronal (somatic, dendritic,
axonal), glial (astrocytic) and extracellular spaces at the $\mu$m
scale (\Cref{tab:scenarios}). These scenarios, along with sample
simulation results, are presented in some detail, with the idea that the
models define and can be reused for future benchmarks. A list of physical
parameters and initial conditions common for all models is provided in
Table~\ref{tab:physconst}. We remark that this setting can be further refined considering heterogeneous initial conditions and physical parameters, depending on the cell type, and heterogeneous membrane properties within a single cell \cite{saetra2021electrodiffusive}. 
All scenarios consider the ionic species $K = \{\mathrm{Na}^+,\mathrm{K}^+, \mathrm{Cl}^-\}$, while $f_r^k=0$ in \eqref{eq:conservation-ions} for all $r$ and $k$ unless otherwise indicated. We denote spatial coordinates by $\xx = (x,y, z)$.
We also refer to the associated
open software repository~\cite{benedusi2024knpemi-zenodo}.
  \begin{table}
    \begin{center}
    \begin{tabular}{lllll}
    \toprule
        Parameter & Symbol & Value & Unit & Ref. \\
    \midrule
        gas constant       & $R$ &  8.314             & J/(K mol)   &  \\
        temperature        & $T$ &  300               & K           & \\
        Faraday constant & $F$ &  $9.648\cdot 10^4$ & C/mol       & \\
        membrane capacitance & $C_M$ & $0.02$ & F& \\
        Na\textsuperscript{+} diffusion coefficient & $D^\text{Na}_r$ & $1.33\cdot10^{-9}$ & m\textsuperscript{2}/s & \cite{hille2001ion}\\
        K\textsuperscript{+} diffusion coefficient  & $D^\text{K}_r$  & $1.96\cdot10^{-9}$ & m\textsuperscript{2}/s & \cite{hille2001ion}\\
        Cl\textsuperscript{\textminus} diffusion coefficient & $D^\text{Cl}_r$ &  $2.03\cdot10^{-9}$ & m\textsuperscript{2}/s & \cite{hille2001ion}\\
        Na\textsuperscript{+} leak conductivity & $g^\text{Na}_\text{leak}$ & 1 & S/m\textsuperscript{2}  & \\
        K\textsuperscript{+} leak conductivity & $g^\text{K}_\text{leak}$  & 4 & S/m\textsuperscript{2} & \\
        Cl\textsuperscript{\textminus} leak conductivity & $g^\text{Cl}_\text{leak}$ & 0   & S/m\textsuperscript{2} & \\
        K\textsuperscript{+}  HH max conductivity & $\bar{g}^\text{K}$ &360  & S/m\textsuperscript{2} & \cite{hodgkin1952quantitative} \\
        Na\textsuperscript{+} HH max conductivity & $\bar{g}^\text{Na}$   & 1200 & S/m\textsuperscript{2} & \cite{hodgkin1952quantitative} \\
        stimulus factor & $\bar{g}_\text{stim}$ & 40 & S/m\textsuperscript{2} & \\
        stimulus time constant & $a$ & 0.002 & s & \\
        initial ICS Na\textsuperscript{+} concentration & $[\text{Na}^+]_i^0$ & 12 & mM & \cite{pods2013electrodiffusion} \\
        initial ECS Na\textsuperscript{+} concentration & $[\text{Na}^+]_e^0$ & 100& mM & \cite{pods2013electrodiffusion} \\
        initial ICS K\textsuperscript{+} concentration  & $[\text{K}^+]_i^0$ & 125 & mM & \cite{pods2013electrodiffusion}\\
        initial ECS K\textsuperscript{+} concentration &$[\text{K}^+]_e^0$ & 4 & mM & \cite{pods2013electrodiffusion}\\
        initial ICS Cl\textsuperscript{\textminus} concentration & $[\text{Cl}^-]_i^0$ & 137 & mM & \cite{pods2013electrodiffusion}\\
        initial ECS Cl\textsuperscript{\textminus} concentration & $[\text{Cl}^-]_e^0$ & 104 & mM & \cite{pods2013electrodiffusion}\\
        initial membrane potential & $\phi_M^0$ & $-67.74$ & mV &  \\
        resting membrane potential & $\phi_{\text{rest}}$ & $-65$ & mV &  \\
        initial Na\textsuperscript{+} activation  & $m_0$ & 0.0379 &  & \cite{hodgkin1952quantitative}\\
        initial Na\textsuperscript{+} inactivation & $h_0$ & 0.688 &  & \cite{hodgkin1952quantitative}\\
        initial K\textsuperscript{+} activation & $n_0$ & 0.276 &  & \cite{hodgkin1952quantitative}\\
        maximum pump rate & $\rho_{\mathrm{pump}}$ &  $1.115\cdot 10^{-6}$ & mol/m$^2$ s & \cite{halnes2013electrodiffusive}\\
        ICS Na\textsuperscript{+} threshold for Na/K-pump & $P_{\mathrm{Na}}$ &  10 & mol/m$^3$  & \cite{halnes2013electrodiffusive}\\
        ECS K\textsuperscript{+} threshold for Na/K-pump & $P_{\mathrm{K}}$ &  1.5 & mol/m$^3$  & \cite{halnes2013electrodiffusive}
    \end{tabular}
    \caption{\label{tab:physconst} Physical parameters and initial values, based on~\cite{ellingsrud2020finite}.}
    \end{center}
\end{table}
\subsection{Model A}
\begin{figure}   
    \begin{subfigure}[b]{0.32\textwidth}
    \centering
    \begin{tikzpicture}[scale=0.8]
    \draw[draw=black,very thick,fill=blue!10!white] (0,0) rectangle (4,4);
    \draw[draw=black,very thick,fill=red!10!white] (1,1) rectangle (3,3);

    \node at (1.7,0.3) {$\mathcal{T}_e$};
    \node at (1.7,2.2) {$\mathcal{T}_i$};
    \node at (0.6,1.5) {$\Gamma$};
    \node at (-0.5,1) {$\partial\Omega$};

    \draw[] (0.8,1.5) -- (1,1.5);
    \draw[] (-0.2,1) -- (0,1);

    \draw[gray, very thin] (0,0) -- (4,4);
    \draw[gray, very thin] (1,0) -- (4,3);
    \draw[gray, very thin] (2,0) -- (4,2);
    \draw[gray, very thin] (3,0) -- (4,1);
    \draw[gray, very thin] (0,1) -- (3,4);
    \draw[gray, very thin] (0,2) -- (2,4);
    \draw[gray, very thin] (0,3) -- (1,4);

    \draw[gray, very thin] (1,0) -- (1,4);
    \draw[gray, very thin] (2,0) -- (2,4);
    \draw[gray, very thin] (3,0) -- (3,4);

    \draw[gray, very thin] (0,1) -- (4,1);
    \draw[gray, very thin] (0,2) -- (4,2);
    \draw[gray, very thin] (0,3) -- (4,3);
    \end{tikzpicture}
    
    \caption{2D geometry}
    \label{fig:model_A_geo_2D}
    \end{subfigure}
    \begin{subfigure}[b]{0.32\textwidth}
      \includegraphics[width=\textwidth]{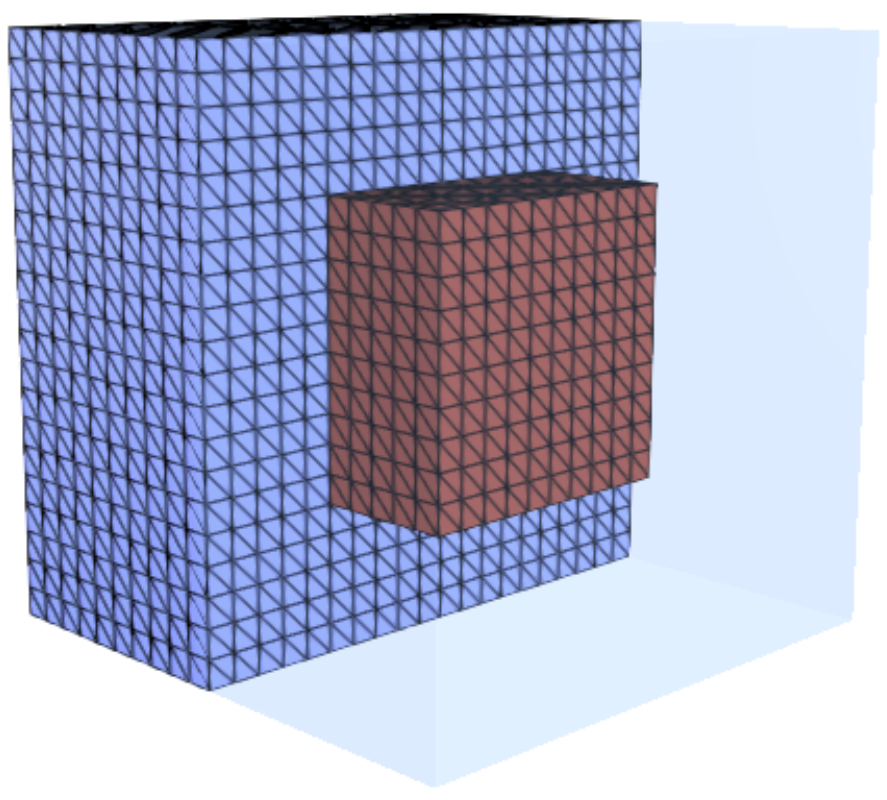}
      \caption{3D geometry}
      \label{fig:model_A_geo_3D}
    \end{subfigure}
    \begin{subfigure}[b]{0.3\textwidth}
    \includegraphics[width=\textwidth]{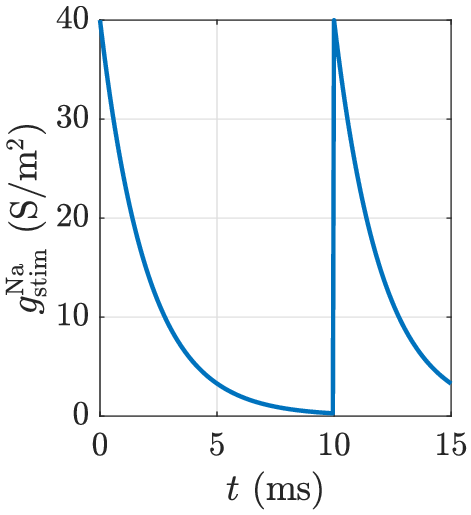}
      \caption{Stimulus}
      \label{fig:model_A_stimulus}
    \end{subfigure}
    \begin{subfigure}[b]{0.32\textwidth}
    \centering
    \includegraphics[width=\textwidth]{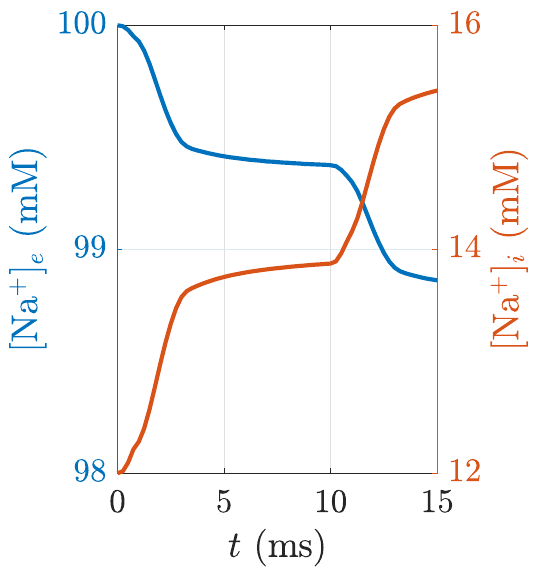}
    \caption{Sodium} 
    \label{fig:model_A_Na}    
    \end{subfigure}
    \hfill
    \begin{subfigure}[b]{0.33\textwidth}
    \centering
    \includegraphics[width=\textwidth]{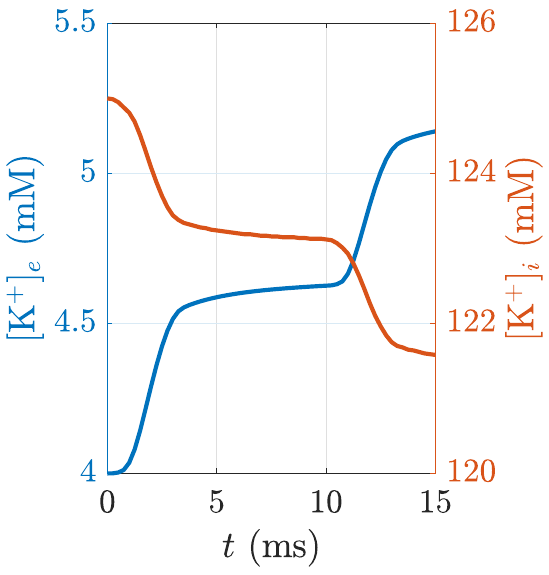}
    \caption{Potassium} 
    \label{fig:model_A_K}    
    \end{subfigure}
    \hfill
    \begin{subfigure}[b]{0.31\textwidth}
    \centering
    \includegraphics[width=\textwidth]{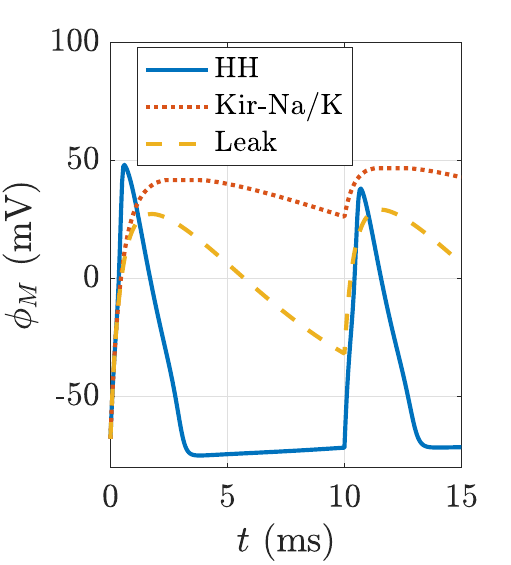}
    \caption{Membrane potential}
    \label{fig:model_A_phi}    
    \end{subfigure}
    \caption{\textbf{Model A}. (a) Sample 2D geometry ($N_x = 4$). (b) Sample 3D geometry ($N_x = 20$). (c) Membrane current stimulus (acting on $\Gamma$) $g_{\rm stim}^{\rm Na}$ versus time $t$. (d) Evolution of ion concentrations over time. Extracellular and intracellular quantities are sampled respectively at $\xx_e=(0.15,0.15)\,\mu\mathrm{m}\in\Omega_e$ and $\xx_i=(0.5,0.5)\,\mu\mathrm{m}\in\Omega_i$. (f) Evolution of membrane potential, which is homogeneous in $\Gamma$, for different ionic models: the default Hodgkin-Huxley (HH); the Kir-Na/K model; a leak model obtained by setting $f\mathrm{_{Kir}}=j_\mathrm{pump}=0$ in the Kir-Na/K model.}
    \label{fig:model_A}
  \end{figure}
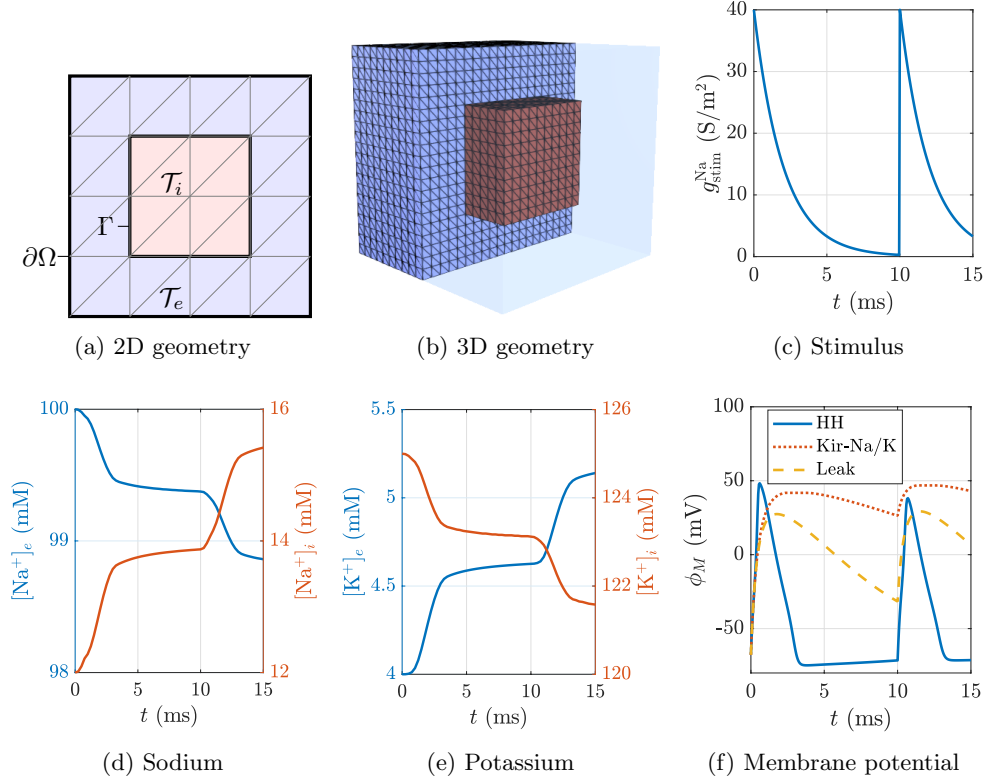

We consider an idealized geometry for a single cell represented by a $d$-dimensional cuboid $\Omega=[0,1]^d$ $\mu$m, with $\Omega_i=[0.25,0.75]^d$ $\mu$m and $d=2,3$ (Figures~\ref{fig:model_A_geo_2D}--\ref{fig:model_A_geo_3D}). We consider the active Hodgkin-Huxley membrane model (Section~\ref{sec:act_dyn}) and, as stimulus, impose a time-periodic sodium current on the entire membrane $\Gamma$:
\begin{equation}
  {g}^{\rm Na}_{\rm stim}={\bar{g}}_{\rm stim}e^{-(t\,\text{mod}\,\tau )/a},
\end{equation}
with $\tau = 10$ ms the time interval between subsequent stimuli (Figure~\ref{fig:model_A_stimulus}). 

We discretize $\Omega$ with a uniform grid with $N_x$ intervals per side and subdivide each cube into 2 triangles in 2D (\Cref{fig:model_A_geo_2D}) or 6 tetrahedra in 3D (\Cref{fig:model_A_geo_3D}). This tessellation results in $N = \left ( ( pN_x+1)^d+o(p^dN_x^d) \right )(|K|+1)$ total degrees of freedom, with $p$ the Lagrange finite element order. The lower order term $o(p^dN_x^d)=(3-d/2)(pN_x)^{d-1} + 2(d-2)$ appears since the degrees of freedom on $\Gamma$ are repeated for $\Omega_i$ and $\Omega_e$. Figures~\ref{fig:model_A_Na}--\ref{fig:model_A_phi} show the
evolution of concentrations and potential over time
for $N_t=300$ time steps, $N_x = 64$, and $p = 1$.

\subsection{Model B}

\begin{figure}
    \centering
    \begin{subfigure}[b]{0.66\textwidth}
    \centering
    \includegraphics[width=\textwidth]{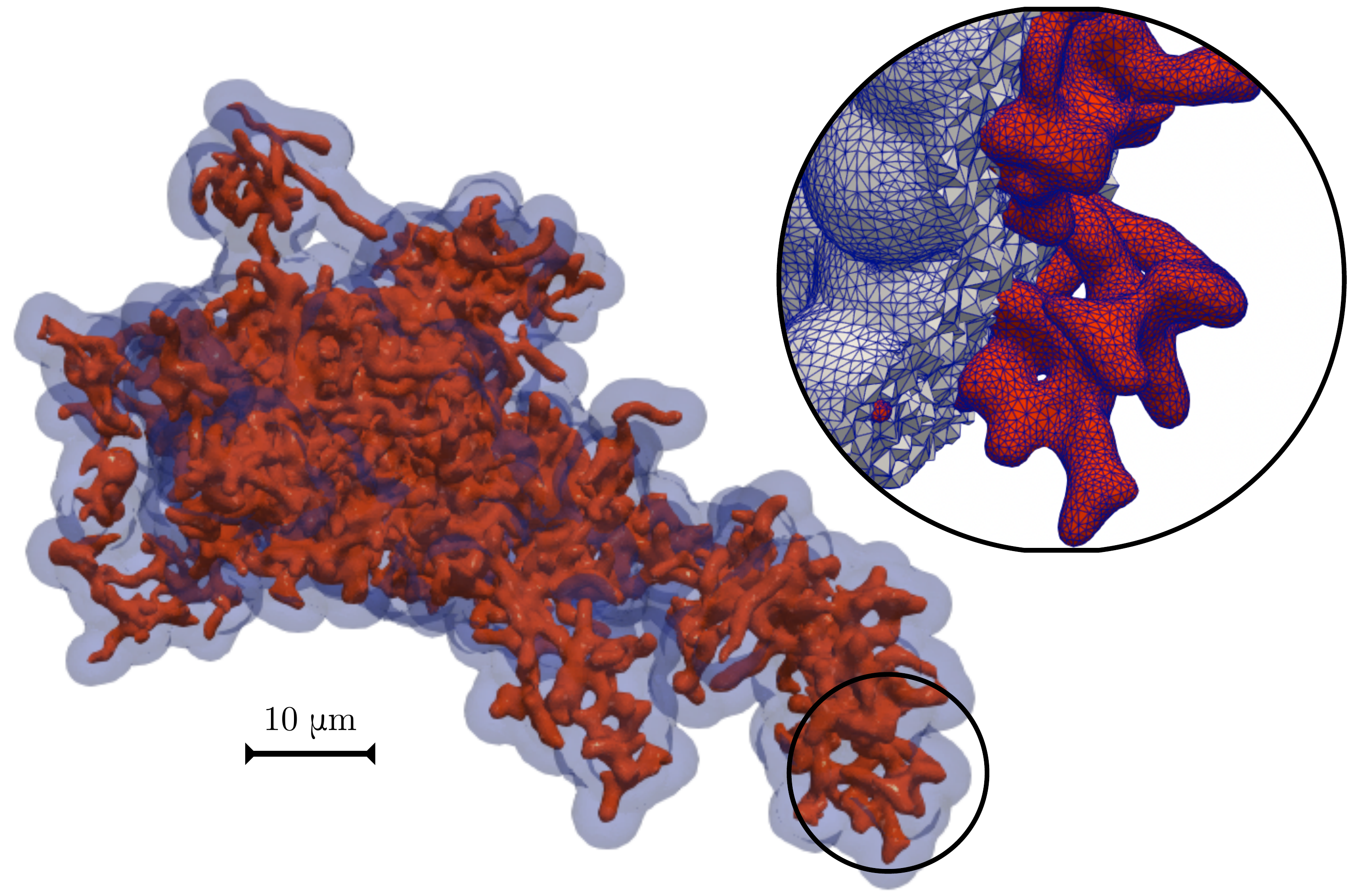}    
    \caption{Astrocyte geometry}
    \label{fig:model_B_geo}
    \end{subfigure}    
    \begin{subfigure}[b]{0.32\textwidth}
    \centering
    \includegraphics[width=\textwidth]{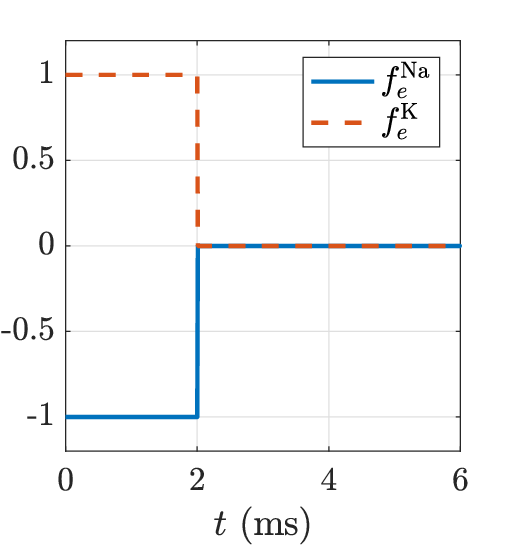}    
    \caption{Stimulus}
    \label{fig:model_B_stim}
    \end{subfigure}    
    \begin{subfigure}[b]{0.44\textwidth}
    \centering
     \includegraphics[width=\textwidth]{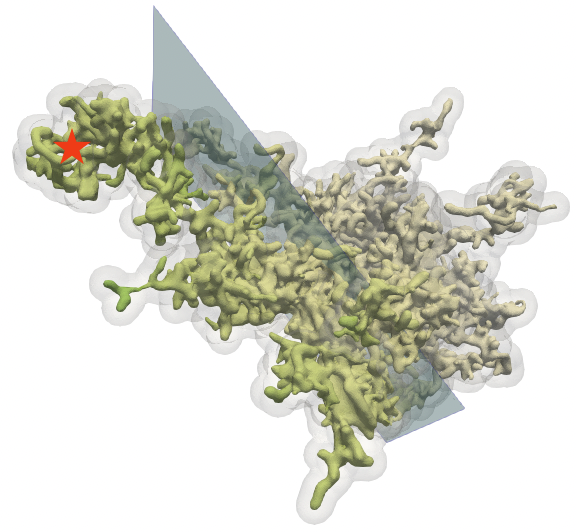}  
    \caption{[K$^+$]$_i$ at $t=30$ s and $x=0$ plane}
    \label{fig:model_B_sol1}
    \end{subfigure}
    \begin{subfigure}[b]{0.55\textwidth}
    \centering
     \includegraphics[width=\textwidth]{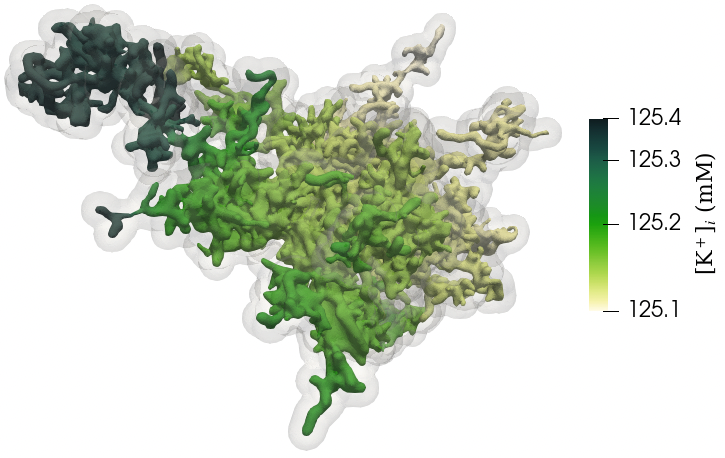}  
    \caption{[K$^+$]$_i$ at $t=200$ ms}
    \label{fig:model_B_sol2}
    \end{subfigure}    
    \begin{subfigure}[b]{0.32\textwidth}
    \centering
    \includegraphics[width=\textwidth]{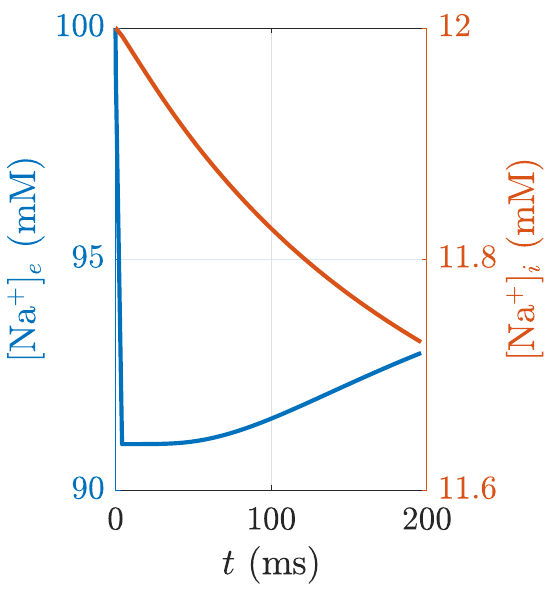}
    \caption{Sodium} 
    \label{fig:model_B_Na}    
    \end{subfigure}
    \begin{subfigure}[b]{0.32\textwidth}
    \centering
    \includegraphics[width=\textwidth]{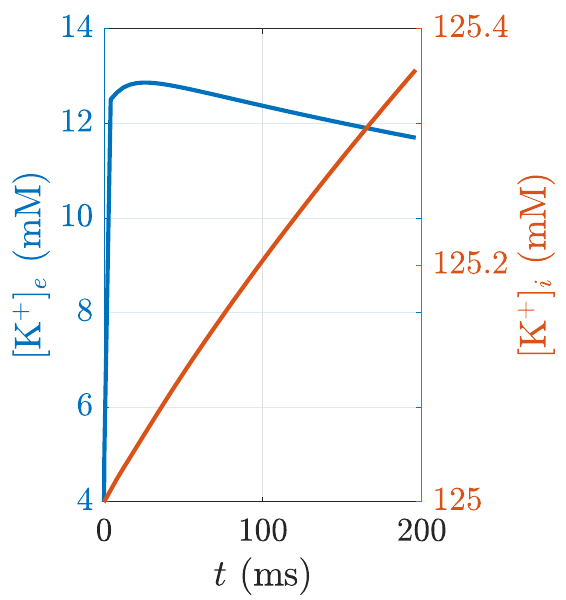}
    \caption{Potassium} 
    \label{fig:model_B_K}    
    \end{subfigure}
    \begin{subfigure}[b]{0.34\textwidth}
    \centering
    \includegraphics[width=\textwidth]{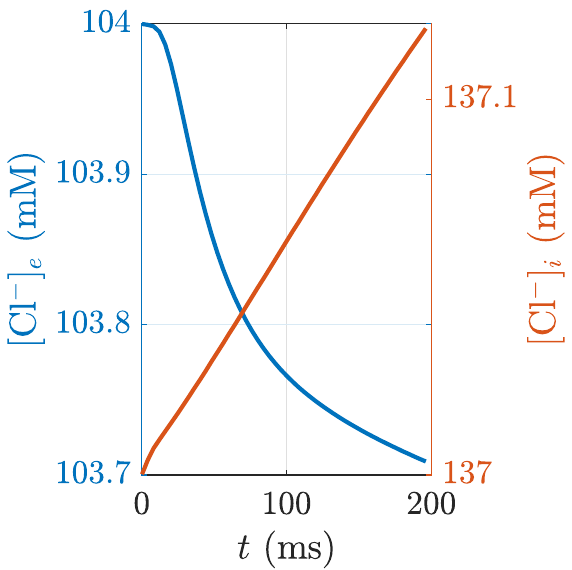}
    \caption{Chloride} 
    \label{fig:model_B_Cl}    
    \end{subfigure}
    \caption{\textbf{Model B}. (a) Astrocyte geometry with corresponding tetrahedral mesh. (b) Ionic source terms over time, imposed for $x<0$. (c),(d) Intracellular potassium [K]$_i$ (mM) at two different time steps. Potassium is injected in the extracellular space for $x<0$. At $t=2$ ms the extracellular potassium [K]$_e$ is increased by approximately 10 mM w.r.t. the initial condition [K]$^0_e=4$ mM. As the system evolves, potassium is electro-diffusing in both $\Omega_i$ and $\Omega_e$ through $\Gamma$. (e)--(g) Time evolution of ionic concentrations $[k]_r(\xx_\star,t)$ (mM) with $\xx_\star\in\Gamma$ labeled by $\star$ in panel (c).}
      \label{fig:model_B}
\end{figure}

We examine the electrodiffusive response of an astrocyte (a
heavily branching-type of glial cell) in response to local changes in
extracellular potassium and sodium concentrations resulting e.g.~from
surrounding neuronal activity. The intracellular astrocytic domain
$\Omega_i$ is described by a 3D
geometry obtained via Ultraliser \cite{abdellah2023ultraliser,marwan_abdellah_2022_7105941}, with a
diameter of approximately 90 $\mu$m, centered at the origin. The extracellular
domain $\Omega_e$ is constructed as a 2 $\mu$m thick layer enclosing
$\Omega_i$ (Figure~\ref{fig:model_B_geo}). The astrocytic membrane
dynamics are governed by the passive Kir--Na/K membrane model
(Section~\ref{sec:pass_dyn}) over the entire interface $\Gamma$. As
stimulus, we consider a bolus injection of potassium and removal of
sodium in a section of the extracellular space over a time period of 2
ms:
\begin{equation*}
f_e^{\rm Na}(\xx, t) = -1,\,f_e^{\rm K}(\xx, t) = 1, \quad \xx \in \Omega_e\quad \text{s.t.}\quad x < 0, \quad 0 < t \leq  2 \,\text{ms} ,
\end{equation*}
and $f_r^k=0$ otherwise,  cf. Figure~\ref{fig:model_B_stim}. No membrane stimulus is applied ($g_{\rm stim}^{\rm Na}=0$).

The geometry is described by a tetrahedral mesh with $N_i = 475 \,
206$ and $N_e = 860 \, 025$ vertices, and a total problem size of $N
\approx 5.3 \cdot 10^6$ for
$p=1$. Figures~\ref{fig:model_B_sol1}--\ref{fig:model_B_Cl} show the
evolution of concentrations over time for $\Delta t=0.1$ ms and $N_t =
2000$.

\subsection{Model C}
\begin{figure}    
    \centering        
    \begin{subfigure}[b]{\textwidth}
    \centering
    \includegraphics[width=0.85
    \textwidth]{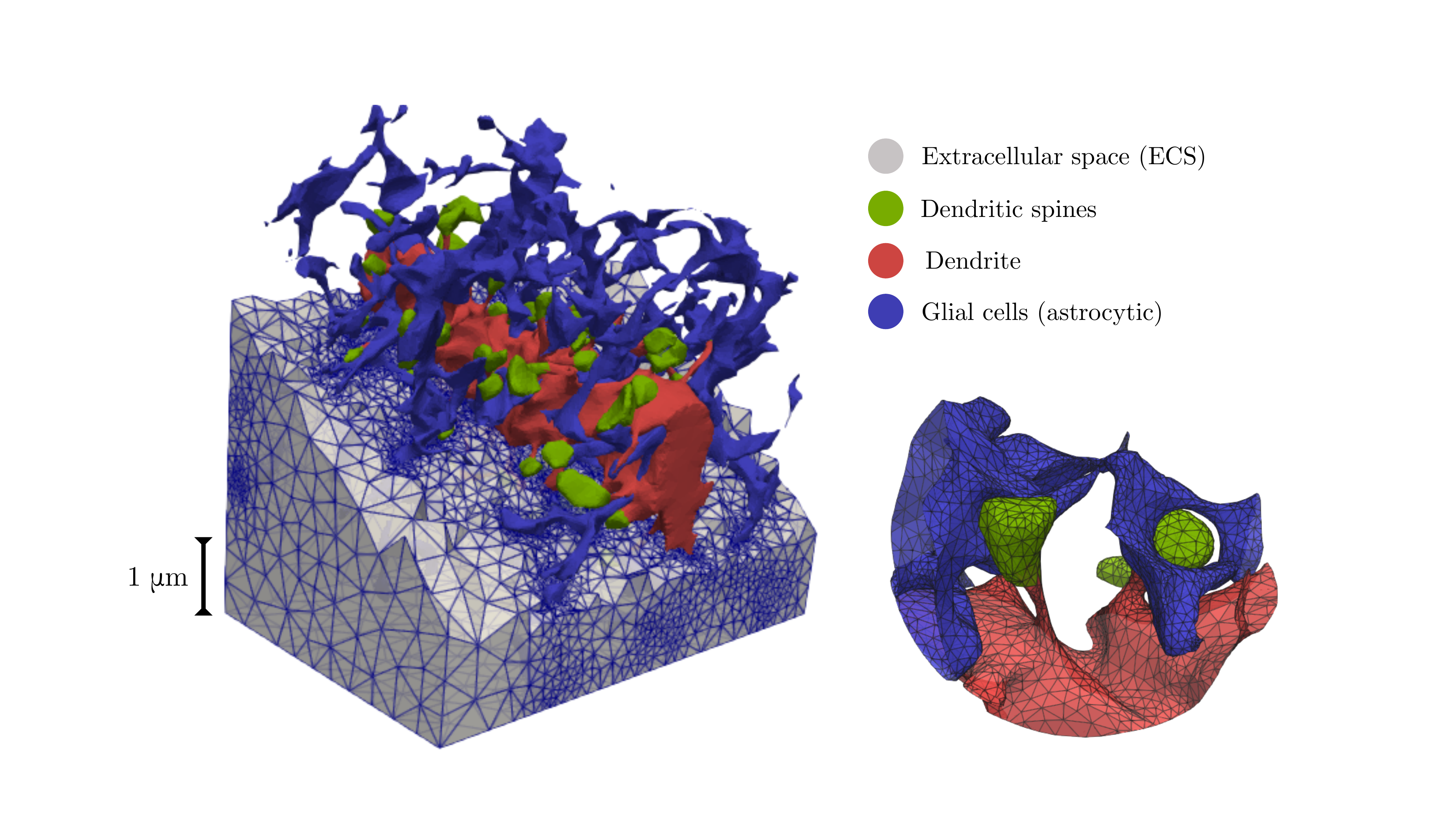}
    \caption{Geometry}
    \label{fig::model_C_geo}
    \end{subfigure}
    \hfill
    \begin{subfigure}[b]{0.43\textwidth}
    \centering
    \includegraphics[width=\textwidth]{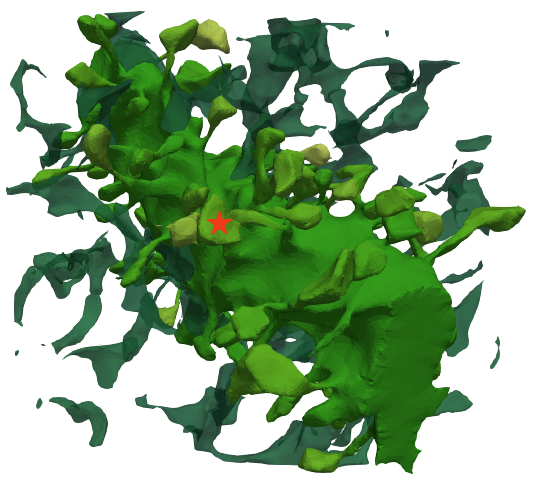}
    \caption{[K$^+$]$_i$ in $\Gamma$ at $t=20$ ms} 
    \label{fig::model_C_sol1}    
    \end{subfigure}
    \begin{subfigure}[b]{0.56\textwidth}
    \centering
    \includegraphics[width=\textwidth]{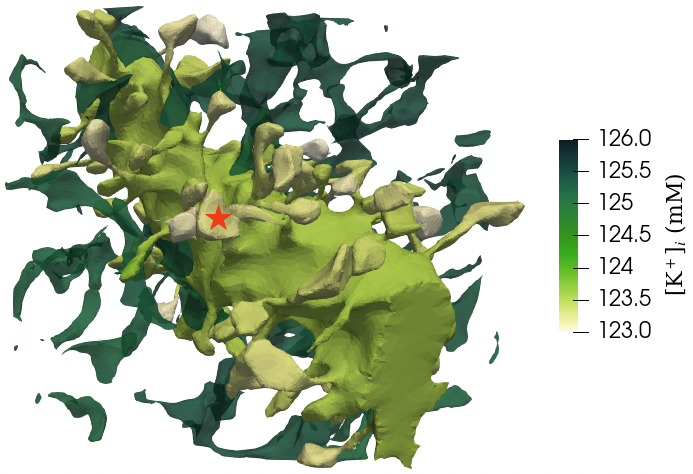}
    \caption{[K$^+$]$_i$ in $\Gamma$ at $t=40$ ms} 
    \label{fig::model_C_sol2}    
    \end{subfigure}
    \begin{subfigure}[b]{0.359\textwidth}
    \centering
    \includegraphics[width=\textwidth]{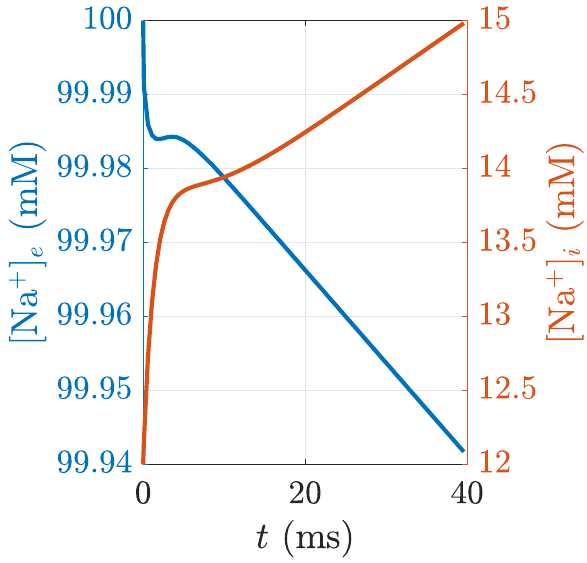}
    \caption{Sodium} 
    \label{fig::model_C_Na}    
    \end{subfigure}
    \begin{subfigure}[b]{0.34\textwidth}
    \centering
    \includegraphics[width=\textwidth]{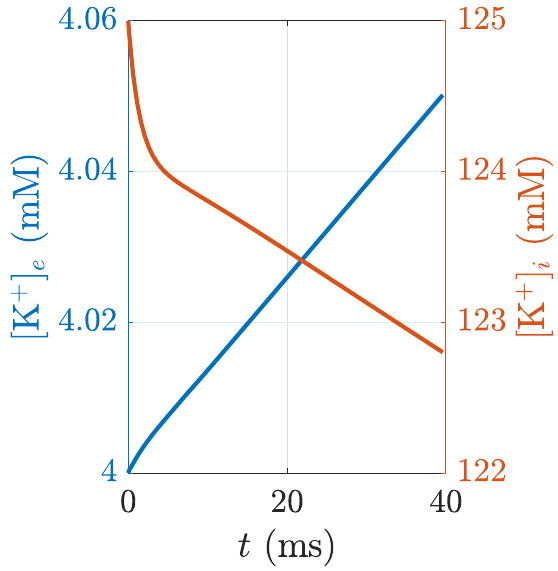}
    \caption{Potassium} 
    \label{fig::model_C_K}    
    \end{subfigure}
    \begin{subfigure}[b]{0.28\textwidth}
    \centering
    \includegraphics[width=\textwidth]{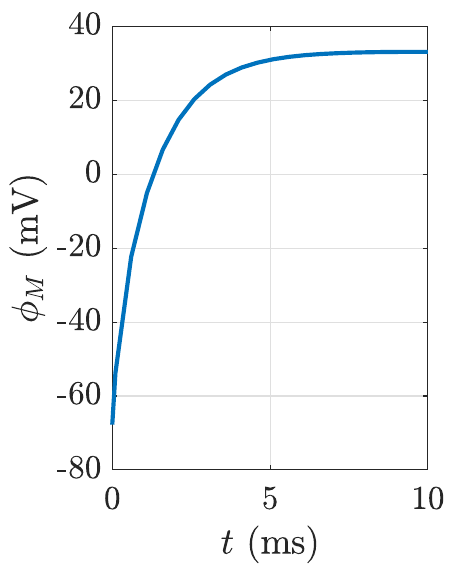}
    \caption{Membrane potential}
    \label{fig::model_C_phi}    
    \end{subfigure}
    \caption{\textbf{Model C}. (a) Dendritic segment (red), spines heads $\Gamma_{\rm head}$ (green), and glial cells (blue). A portion of the meshed extracellular domain is shown as long as a particular of few meshed dendritic spines. (b),(c) [K]$_i$ (mM) for two different time steps, given a membrane stimulus located in the spines heads. (d)--(f) Concentrations $[k]_r(\xx_\star,t)$ (mM) for $k\in\{\text{K$^+$, Na$^+$}\}$ and potential $\phi_M(\xx_\star,t)$ (mV) with $\xx_\star\in\Gamma_{\rm head}$. The location of $\xx_\star\in\Gamma_{\rm head}$ is labeled in panels (b) and (c) with $\star$.}
    \label{fig::model_C}
\end{figure}

This model focuses on the electrodiffusive interplay between neuronal
structures, the extracellular space and surrounding glial
structures at the scale of dendritic spines. We consider a 3D geometry
(\Cref{fig::model_C_geo}) representing a neuronal dendrite segment
with multiple dendritic spines, surrounded by extracellular space and
glial cells, all immersed in a
$8\times8\times5$ $\mu$m cuboid~\cite{kinney2013extracellular,
  bartol2015computational}. We remark that even this highly
non-trivial geometry represents a substantial simplification of the
dense tissue structure, and in particular that a large number of
glial structures are ignored. We apply the passive Kir--Na/K dynamics
model (Section~\ref{sec:pass_dyn}) at both the astrocyte--ECS and
dendrite--ECS interfaces, though with $f_{\mathrm{Kir}} =
j_{\mathrm{pump}} = 0$ for the latter, thus in essence only
considering leak channel dynamics on the dendrite--ECS interfaces. To trigger a system response,
we apply sodium membrane currents at the dendritic spine heads
$\Gamma_{\rm head} \subset \Gamma$ (marked in green in \Cref{fig::model_C_geo}), mimicking synaptic signalling:
\begin{equation}
  {g}^{\rm Na}_{\rm stim}(\xx)=
  \begin{cases}2\bar{g}_{\rm stim}~&{\text{ if }}~\xx\in\Gamma_{\rm head},\\0~&{\text{ if }}~\xx\notin \Gamma_{\rm head},\end{cases}
\end{equation}
We consider a tetrahedral mesh with $N_i=193\,129$ and $N_e=306\,478$ vertices, and a total problem size of $N\approx2.0\cdot10^6$ for $p=1$. In Figures~\ref{fig::model_C_sol1}--\ref{fig::model_C_phi}, we show the resulting concentrations and potential for $N_t = 400$ time steps and $\Delta t = 0.1$ ms.

\subsection{Model D}

\begin{figure}
    \centering       
    \begin{subfigure}[b]{\textwidth}
    \centering
    \includegraphics[width=0.8\textwidth]{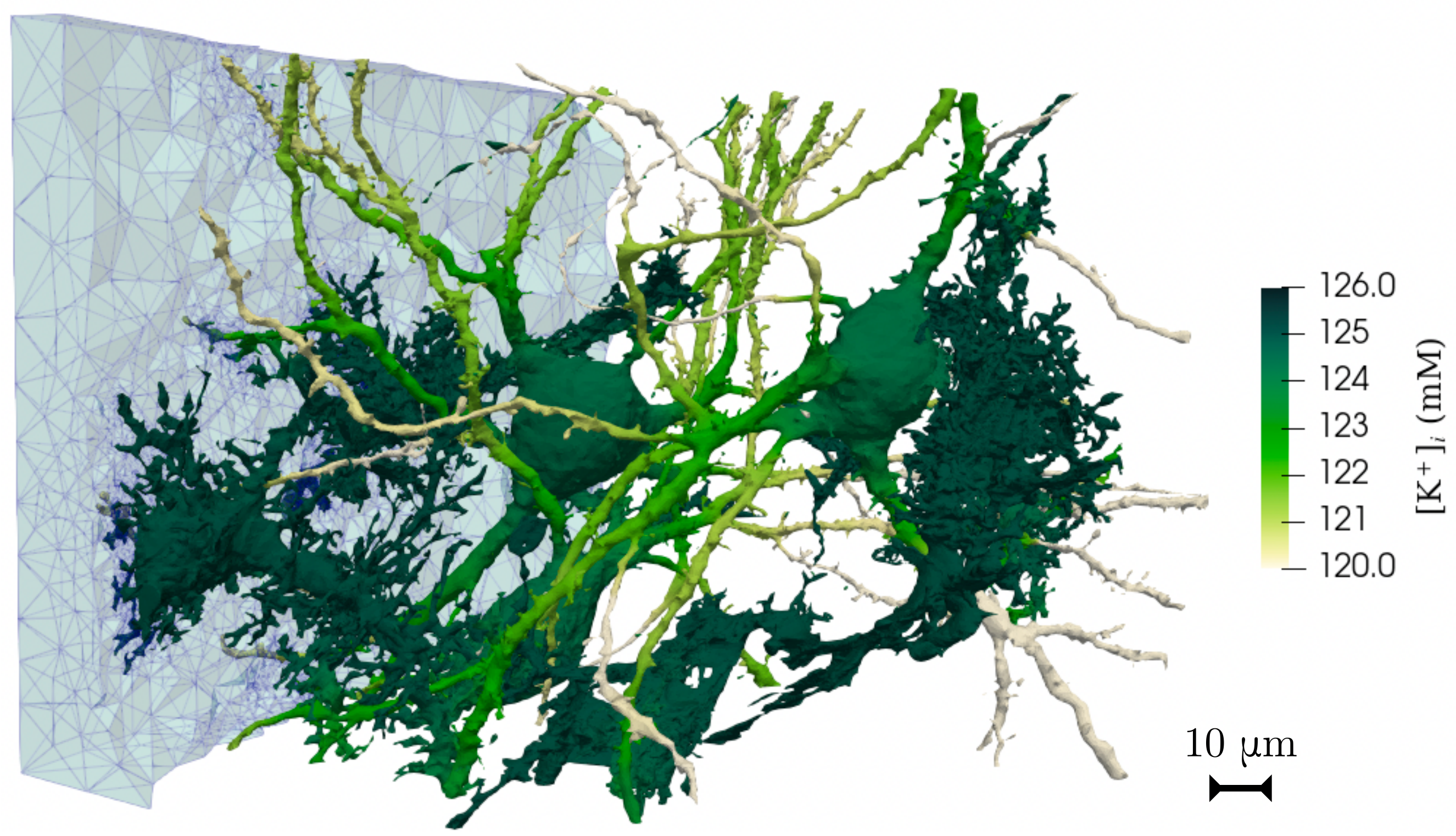}    
    \caption{Geometry and intracellular potassium in $\Gamma$ at $t=50$ ms}
    \label{fig:model_D_geo}
    \end{subfigure}
    \hfill    
    \begin{subfigure}[b]{0.43\textwidth}
    \includegraphics[width=\textwidth]{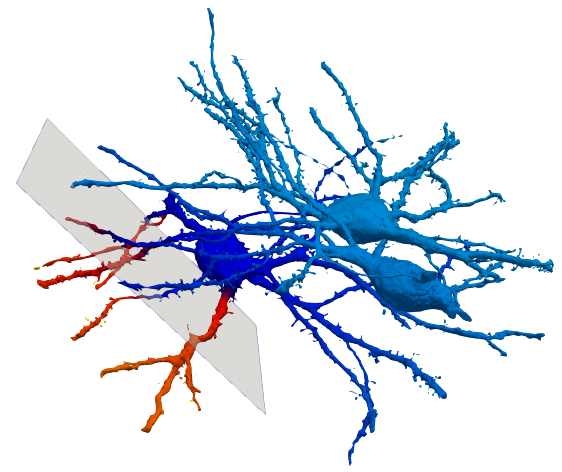}    
    \caption{Membrane potential at $t=0.5$ ms}
    \label{fig:model_D_sol1}    
    \end{subfigure}            
    \begin{subfigure}[b]{0.555\textwidth}    
    \includegraphics[width=\textwidth]{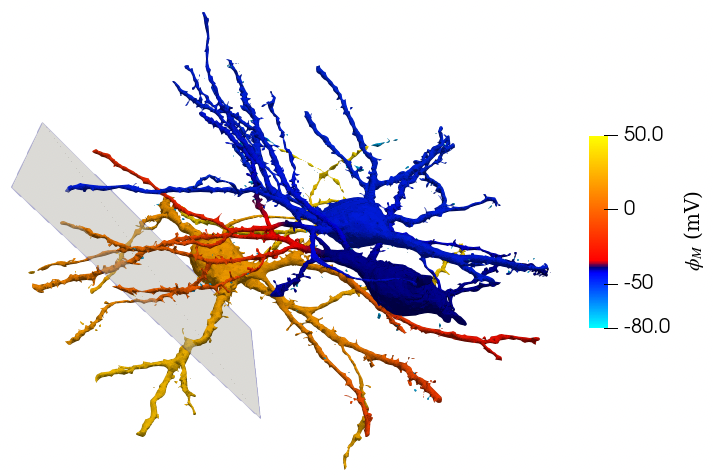}    
    \caption{Membrane potential at $t=1$ ms}
    \label{fig:model_D_sol2}        
    \end{subfigure}            
    \begin{subfigure}[b]{0.328\textwidth}
    \centering
    \includegraphics[width=\textwidth]{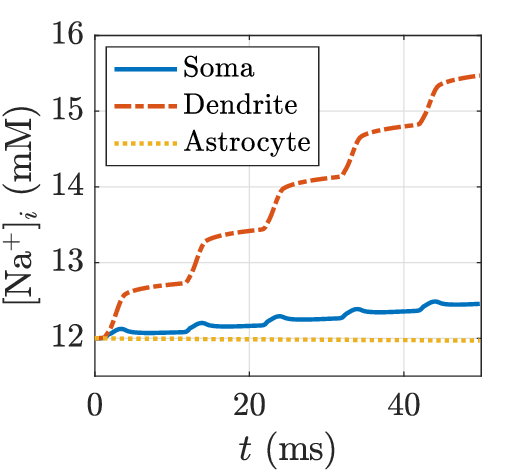}
    \caption{Sodium} 
    \label{fig::model_D_Na}    
    \end{subfigure}
    \begin{subfigure}[b]{0.328\textwidth}
    \centering
    \includegraphics[width=\textwidth]{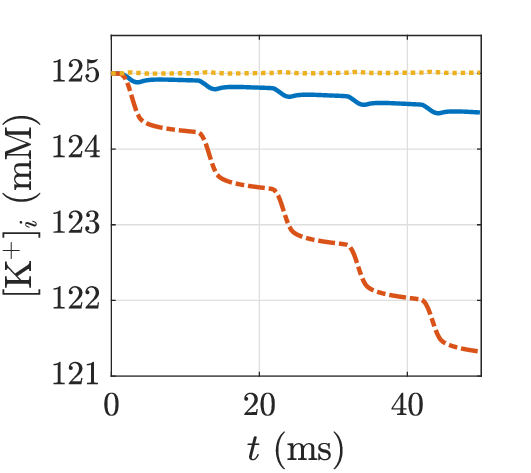}
    \caption{Potassium} 
    \label{fig::model_D_K}    
    \end{subfigure}
    \begin{subfigure}[b]{0.328\textwidth}
    \centering
    \includegraphics[width=\textwidth]{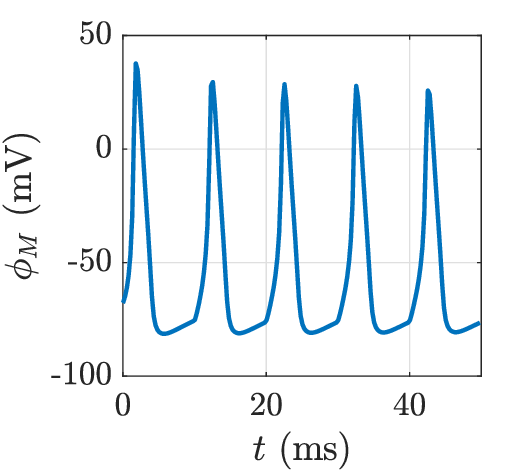}
    \caption{Membrane potential}
    \label{fig::model_D_phi}    
    \end{subfigure}
    \caption{\textbf{Model D}. (a) Cortex geometry with neurons and astrocytes, with a portion of meshed ECS. Intracellular potassium at $t=50$ ms is shown. (b),(c) Spreading electric potential on neurons from an initial localized stimulus for two different time steps. A portion of the plane $x = 10\,\mu$m is shown, as boundary of the stimulus region. A threshold potential of 40 mV is highlighted by the color scale transition. (d),(e) Intracellular concentrations for $N_t=1000$. The solution is sampled in the membrane of a dendrite, a soma, and an astrocyte. Legend is shared between the two images. (f) Evolution of neuronal membrane potential.}
    \label{fig:model_D}
\end{figure}

We examine electrodiffusive neuronal--ECS--astrocyte interplay at a larger spatial scale. From a $100\times100\times80 \, \mu$m portion of a rat somatosensory cortex \cite{marwan_abdellah_2022_7105941,cali20193d}, we select three neurons and four astrocytes and embed these in a cuboid extracellular space (\Cref{fig:model_D_geo}). For the neuronal membranes, we apply the active Hodgkin-Huxley model (Section~\ref{sec:act_dyn}) with a periodic, localized stimulus current
\begin{equation}
  g^{\rm Na}_{\rm stim}=10\cdot{\bar{g}}_{\rm stim}e^{-(t\,\text{mod}\,\tau )/a} \quad \text{for} \quad x < 10\,\mu\text{m},
\end{equation}
and ${g}^{\rm Na}_{\rm stim}=0$ otherwise. The period $\tau = 10$ ms. For the astrocytic membranes, we use the passive Kir--Na/K-model as previously and ${g}^{\rm Na}_{\rm stim}=0$ (Section~\ref{sec:pass_dyn}). 

We consider tetrahedral meshes with up to $N_i=9\,390\,510$ and $N_e=15\,541\,672$ and problem size $N\approx1.0\cdot10^8$. The numerical solution is shown in Figures~\ref{fig:model_D_sol1}-\ref{fig::model_D_phi}. 

\section{Numerical experiments}
\label{sec::exp}

In this section, we examine the robustness and efficiency of the
proposed monolithic numerical strategy in terms of iteration count and
parallel performance when applied to Models A--C. All numerical
experiments are performed on the Norwegian supercomputer
Saga\footnote{\url{https://documentation.sigma2.no/hpc_machines/saga.html}},
using up to 256 Intel Xeon-Gold 6138/6230R 2.0-2.1 GHz CPU cores with
40 cores per node. 

\subsection{Implementation and numerical verification}

Our KNP-EMI solver implementation~\cite{benedusi2024knpemi-zenodo} is
based on the well-established FEniCS finite element software~\cite{alnaes2015fenics,logg2012automated}, using the multiphenics library~\cite{multiphenics} to handle variational problems coupled across multiple subdomains. The linear algebra backend is PETSc~\cite{petsc-web-page}. The implementation is flexible in terms of ionic models of arbitrary complexity, and has been verified using the method of manufactured solutions for the same benchmark problems as introduced in~\cite[Section 3.1]{ellingsrud2020finite}, obtaining the expected convergence rates. The image-based volumetric meshes (in Models B--D) are generated from segmented microscopy data via unions of lower-dimensional surfaces using the computational geometry library fTetWild \cite{ftet}. 

\subsection{Solver configurations}
\label{sec:imple}

When using a GMRES method to solve the KNP-EMI system, we restart the algorithm after 30 iterations. For each time step $n$, we use the solution at $n-1$ as an initial guess, substantially reducing the number of GMRES iterations. As stopping criterion for the iterative solvers, we consider a tolerance of $10^{-6}$ for the preconditioned relative residual. To evaluate $P_0^{-1}$, we use a single V-cycle of hypre BoomerAMG \cite{falgout2002hypre}. We refer to the associated software repository~\cite{benedusi2024knpemi-zenodo} for the exhaustive list of AMG parameters\footnote{For 3D runs, we set the threshold for strong coupling to 0.5, instead of the default value (0.25).}. In terms of computational time to solve the system, the GMRES method is preferable compared to BiCGStab for all the reported experiments. For comparative tests, we also use MUMPS as a direct solver.

To solve the gating variables ODEs (cf.~\eqref{eq::gating}), we use $N_{\text{ode}}=25$ intermediate time steps of the Rush-Larsen method for each time interval, set to $\Delta t = 0.05$ ms (unless noted otherwise). In terms of runtime, solving the ODEs is negligible compared to the linear system solve.

\subsection{Iterative solver robustness}
\label{sec::results}

We begin by studying the robustness of the proposed preconditioning strategy with respect to the spatial and temporal discretization. We consider a single time step for Model A, with $d=2$, as a flexible benchmark example, and vary the time step $\Delta t$, the finite element (polynomial) degree $p\in\{1, 2\}$ and the mesh resolution $N_x \in \{16, 32, 64, 128, 256, 512\}$. The resulting linear system sizes vary from $N=4.6\cdot10^3$ for $(N_x,p)=(16,1)$ to $N=4.2\cdot10^6$ for $(N_x,p)=(512,2)$. For $\Delta t \in \{1, 10, 100\}$ (ms), GMRES preconditioned by $P_0$ (LU($P_0$)) converges in 3 iterations for all cases. For $\Delta t \in \{0.01, 0.1 \}$ (ms), each case converges in 4 iterations. In contrast, the unpreconditioned GMRES method never converges in less than 1000 iterations.

Next, we investigate how to approximate $P_0^{-1}$ efficiently. For larger problems, it is not convenient to invert $P_0$ to full precision. Instead, a suitable multilevel approximation of $P_0^{-1}$ can reduce the time-to-solution while retaining the rapid convergence of the GMRES method. We here compare a set of such approximation strategies to black box approaches and the exact preconditioners for increasing problem sizes $N$, cf. Table~\ref{tab:mod_A_runtimes}. 
We immediately observe that the direct solver and ILU(0) or AMG preconditioners are not viable for larger $N$. We therefore consider two approaches where $P_0^{-1}$ is evaluated exactly via LU or CG (LU($P_0$) and CG($P_0$), respectively). Here, the CG solver is applied block-wise in a field-split fashion with $2(|K|+1)$ fields, and preconditioned with AMG. In both cases, the convergence behavior is dramatically improved: the number of iterations no longer increase with the problem size and remains low ($\approx4$). However, we observe only marginal gains in terms of runtime compared to ILU(0). Finally, we consider strategies based on the multilevel approximation of $P_0$, obtained with a single AMG V-cycle (AMG$(P_0)$). Besides a monolithic approach, we consider either block-wise or field-split application of AMG (AMG$_{\text{FS}}(P_0)$). We observe that both of these latter strategies maintain the robustness in terms of a constant and low iteration count ($\approx4$) for all problem sizes. Moreover, we observe a significant improvement in run times for larger problem sizes. The best performance is observed for AMG applied block-wise, with a $70\%$ reduction in runtime compared to the exact preconditioner (LU($P_0$)). 
\begin{table}
    \centering
    \begin{tabular}{l|rrrr}
      \toprule
      $N$ & 17\,412 & 67\,588 & 266\,244 & 1\,056\,772\\
      ($N_x$) & (64) & (128) & (256) & (512) \\
      \midrule
      Assembly  & 2.5 & 8.2  & 31.3 & 123.3 \\         
      \midrule
      Direct (MUMPS)   & 2.1 & 10.0 & 49.5 & 259.1 \\      
      AMG        & 2.3 [22.5] & 24.8 [61.1] & 252.5 [157] & 3697 [585]\\
      ILU(0)     & 0.4 [15.6] & 2.4 [29.5] & 17.1 [63.2] & 176.7 [169]\\              \midrule
      CG($P_0$)  & 3.0 [4.0] & 10.5 [4.0] & 39.8 [4.0] & 174.4 [4.0] \\    
      LU($P_0$)  & 0.5 [4.0] & 3.5 [4.0] & 18.1 [4.0] & 122.4 [4.0] \\    
      AMG($P_0$)   & 0.8 [4.3] & 3.3 [4.1] & 13.8 [4.0] & 55.2 [4.0] \\
      AMG$_{\text{FS}}(P_0)$   & 0.6 [4.0] & 2.2 [4.0] & 9.0 [4.0] & 37.0 [4.0]
    \end{tabular}
    \caption{runtimes (in serial) and when applicable, in square brackets, average number of preconditioned GMRES iterations to convergence over $n = 1, \ldots, N_t$, for preconditioners (Model A, $p = 1$, $d=2$, $N_t=10$ time steps). In the last two rows different inexact preconditioners are used.}
    \label{tab:mod_A_runtimes}
\end{table}

\begin{figure}[H]
    \centering
    \includegraphics[width=0.95\textwidth]{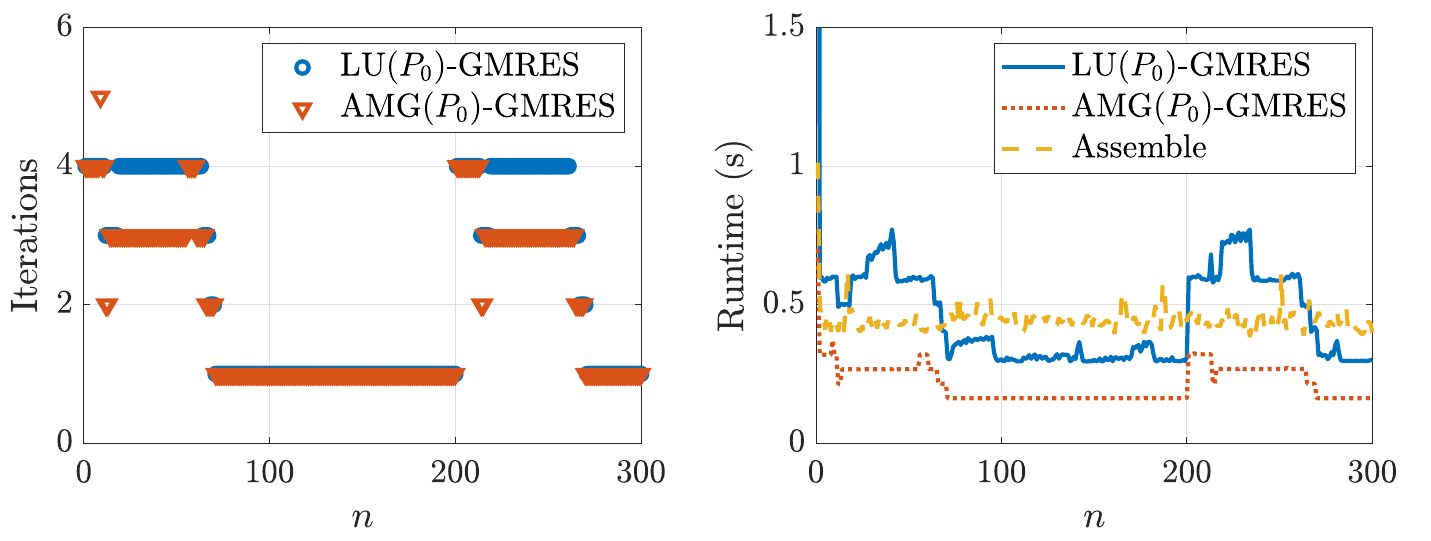}
    \caption{Number of iterations (left) and runtimes (right) over
      $N_t=300$ time steps for Model A with $d = 2$ and $(N_x, p) = (512, 1)$ ($N=1\,056\,772$), using 32 MPI processes.}
    \label{fig:model_A_conv2}
\end{figure}

We end this section by examining how the iteration counts and runtimes depend on time evolution, focusing on the two solution strategies (LU($P_0$) and AMG($P_0$)), now distributed over 32 cores (Figure~\ref{fig:model_A_conv2}). For both strategies, the number of iterations remains low (1--6 iterations) over the whole time span simulated. Moreover, the run time per time step (after initialization) remains bounded throughout the time course, with AMG($P_0$) being approximately $50\%$ faster than LU($P_0$) on average. 

\subsection{Parallel scalability: strong scaling}
\begin{figure}
    \centering  
    \begin{subfigure}[b]{0.495\textwidth}
    \centering
    \includegraphics[width=\textwidth]{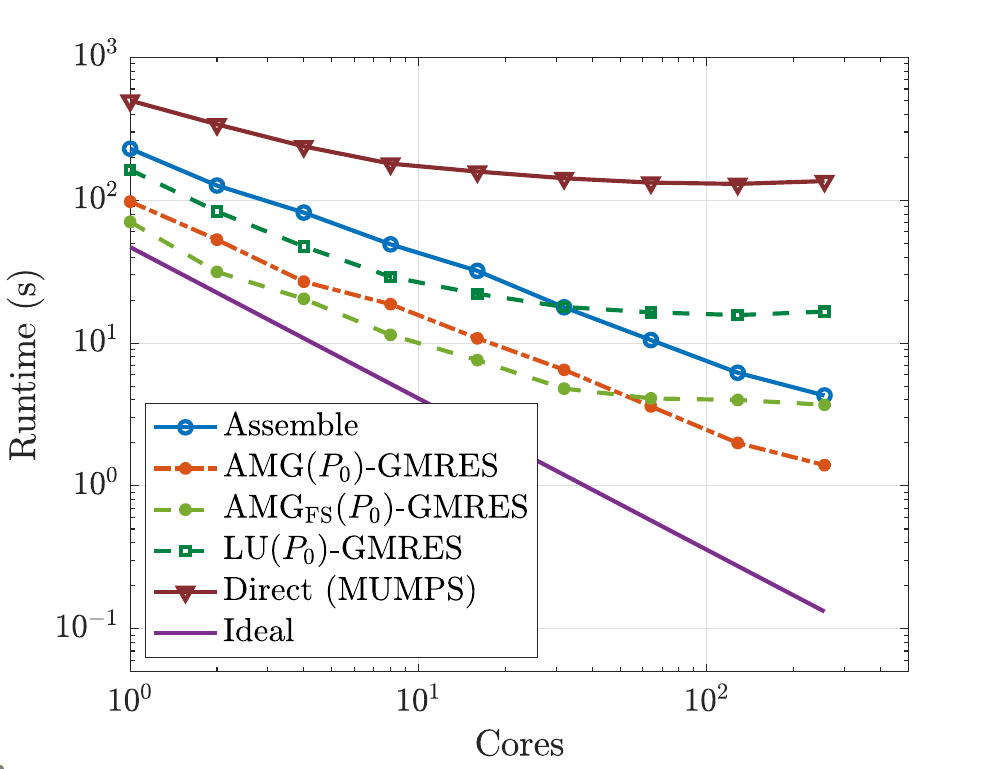}
    \caption{Model A ($d = 2$), $N=1.0\cdot 10^6$ }
    \label{fig::strong_scaling_square}
    \end{subfigure}   
    \hfill
    \begin{subfigure}[b]{0.495\textwidth}
    \centering
    \includegraphics[width=\textwidth]{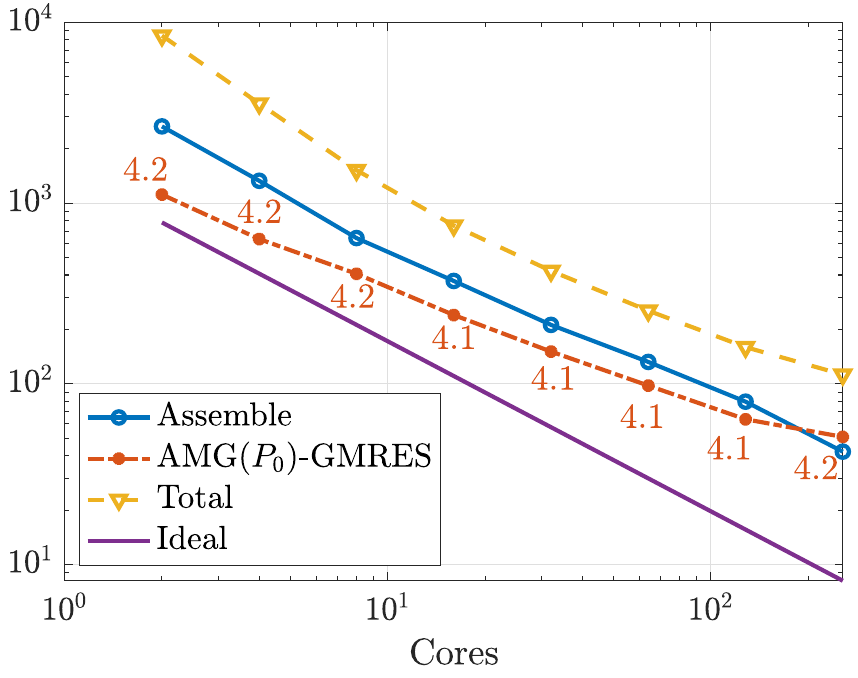}
    \caption{Model A ($d=3$), $N=4.2\cdot 10^6$}
    \label{fig::strong_scaling_cube}
    \end{subfigure}        
    \begin{subfigure}[b]{0.50\textwidth}
    \centering
    \includegraphics[width=\textwidth]{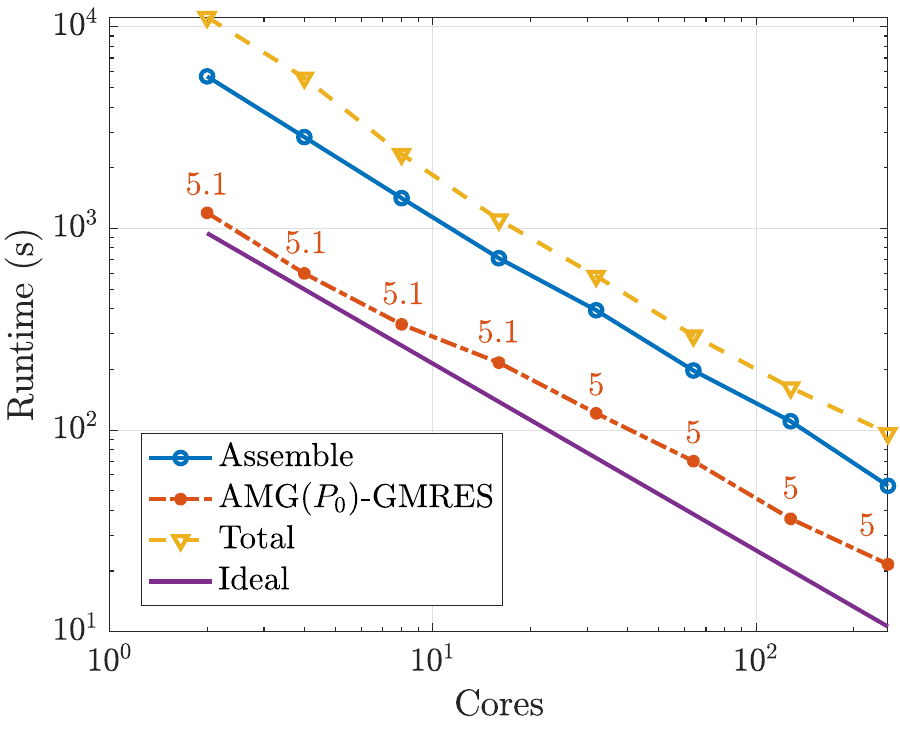}
    \caption{Model B, $N = 5.3 \cdot 10^6$ }
    \label{fig::strong_scaling_B}
    \end{subfigure}    
    \hfill
    \begin{subfigure}[b]{0.48\textwidth}
    \centering
    \includegraphics[width=\textwidth]{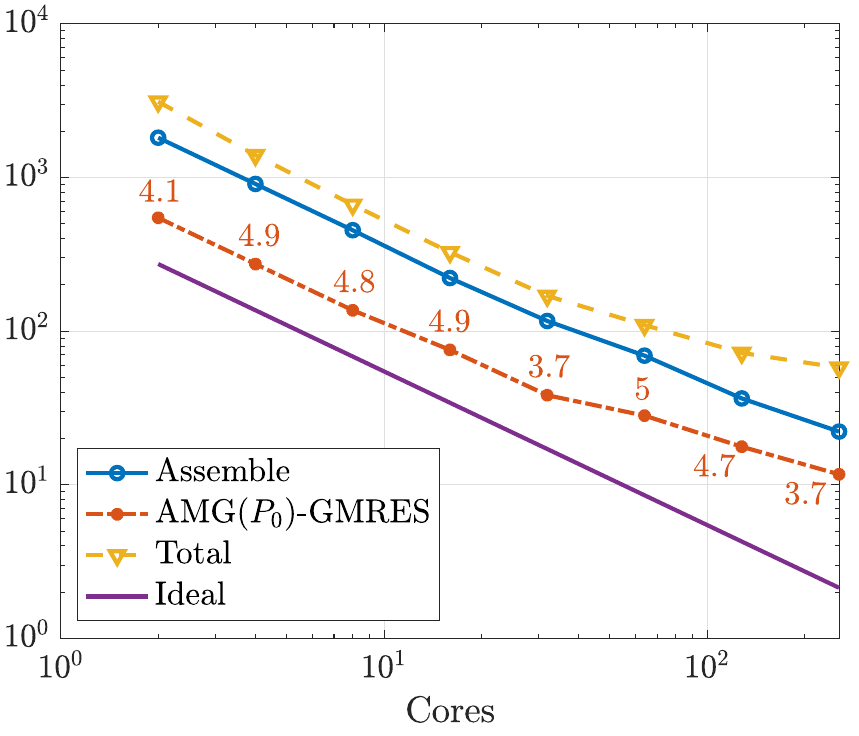}
    \caption{Model C, $N = 2.0 \cdot 10^6$}
    \label{fig::strong_scaling_C}
    \end{subfigure}    
    \caption{Parallel performance: strong scaling for Models A, B, and C. Plots show runtimes for iterative solves, finite element assembly and total runtime versus number of MPI processes (cores). For all models, we use $p = 1$ and $N_t = 20$ time steps. Average AMG($P_0$)-GMRES iterations over $n=1, \ldots, N_t$ are indicated for each data point for Models A ($d=3$), B, and C.}    
    \label{fig::strong_scaling}
\end{figure}
We now more carefully investigate the parallel performance of the more promising preconditioned GMRES strategies (Figure~\ref{fig::strong_scaling}). We first compare their strong scaling (increasing the number of cores while keeping the problem size fixed) for Models A, B, and C. For the smaller problem (Model A with $d=2$), we compare with the direct solver for reference. Throughout, we also examine the scalability of the finite element assembly. Initial tests (using Model A in 2D) indicate that all of the iterative strategies tested scale reasonably well for a few (4--8) cores (Figure~\ref{fig::strong_scaling_square}). The field-split AMG$_{\text{FS}}(P_0)$ is the most efficient for up to $32$ cores, but stagnates for larger core counts since each diagonal block is solved in parallel. The monolithic AMG strategy scales better and outperforms all other strategies for the Model A, $d=2$ test case for $64$ cores or more. We therefore focus our attention on the monolithic AMG in the subsequent experiments. For the larger problems (Models A with $d=3$, B and C), the scaling of AMG($P_0$), the assembly algorithm, and the total runtime are close to ideal for at least up to 256 cores 
(Figures~\ref{fig::strong_scaling_cube}--\ref{fig::strong_scaling_C}). 
We also note that the average number of iterations remain in the same low range (4--5) varying the number of cores. 

\begin{table}
    \centering    
    \begin{tabular}{l|rrrr}
      \toprule
       $N$ & $6.7\cdot10^6$ & $4.5\cdot10^7$ & $4.5\cdot10^7$ & $1.0\cdot10^8$ \\
       $N_i + N_e$  & $1.7\cdot10^6$ & $1.7\cdot10^6$ & $1.3 \cdot 10^7$ & $2.5\cdot10^7$ \\
       Finite element degree ($p$) & $1$ & $2$ & $1$ & $1$ \\
      \midrule
      Initialization ($128$ cores)  & 20.2      &     60.5          & 335.1         & 991.1              \\      
      Initialization ($256$ cores)  & 19.0      & 31.5          & 143.1         & 317.1        \\ 
      \midrule
      Assembly ($128$ cores)        & 82.5      &   511.5            & 468.4         & 1064.1              \\             
      Assembly ($256$ cores)        & 44.0      & 278.5         & 263.2         & 547.2         \\        
      \midrule
      AMG($P_0$) ($128$ cores)      & 20.8 [4.0]& 211.6 [4.0]              & 132.4 [4.6]   & 441.8 [8.0]              \\          
      AMG($P_0$) ($256$ cores)      & 12.0 [3.9]& 111.5 [3.9]   & 64.9 [3.9]    & 242.9 [7.8]  \\   
      \bottomrule
    \end{tabular}
      \caption{Scalability and performance for Model D. Computational runtimes and, in square brackets, average number of iterations to convergence over $n = 1, \ldots, N_t=10$ time steps. Total number of degrees of freedom $N$, intracellular and extracellular degrees of freedom $N_i + N_e$, finite element (polynomial) degree $p$, with rows corresponding to 128 and 256 MPI processes as labeled. 
      The initialization time corresponds to the finite element discretization setup in the first time step, thus becoming negligible for large $N_t$. The columns give (from left to right) the results corresponding to Model D with: (i) $p = 1$, (ii) $p = 2$, (iii) after uniform mesh refinement, and (iv) after uniform and then interface-based mesh refinement.}
      \label{tab:mod_D_runtimes}
\end{table}
We also conduct a focused larger-scale experiment using Model D with up to $10^8$ degrees of freedom to study the effect of bulk- versus interface-dominated mesh refinement on algorithmic and parallel performance (Table~\ref{tab:mod_D_runtimes}). As a baseline, we consider Model D with the default configuration (Table~\ref{tab:mod_D_runtimes}, first column), as described in Section~\ref{sec::models}, for 128 and 256 cores. 
For the case $p = 2$ (Table~\ref{tab:mod_D_runtimes}, second column) with a $6.8\times$ increase in the number of degrees of freedom and reduced sparsity, assembly time increases by a factor of $6.3\times$, while the linear solver time increases by a factor of $9.3\times$. 
On the other hand, after uniform mesh refinement (third column), which yields the same increase in $N$ and same sparsity pattern, the assembly time increases by a factor of $5.7\times$ (128 cores) and $6.0\times$ (256 cores), while the GMRES runtime increases by $6.3\times$ and $5.4\times$, demonstrating near-optimal or super-optimal scalability with respect to $N$.
For both cases, we again observe that the average number of GMRES iterations stays near constant, remaining at 4--5.

Finally, we study the effect of refining the mesh near the interface $\Gamma$ only. This scenario challenges the spectral analysis~\cite{benedusi2021fast} in which the robustness of the preconditioning strategy eliminating membrane blocks is established under the assumption of uniform refinement. Indeed, the average number of iterations now increases, from 4--5 to 8, with a corresponding increase in computational cost (comparing the third and fourth columns of Table~\ref{tab:mod_D_runtimes}).  For all four columns, both the finite element assembly and the AMG($P_0$)-GMRES solver scale near-optimally when increasing from 128 to 256 cores, with an average parallel efficiency of $93\%$ and $94\%$ respectively.

\subsection{Parallel scalability: weak scaling}
\begin{figure}
    \centering  
    \begin{subfigure}[b]{0.51\textwidth}
    \centering
    \includegraphics[width=\textwidth]{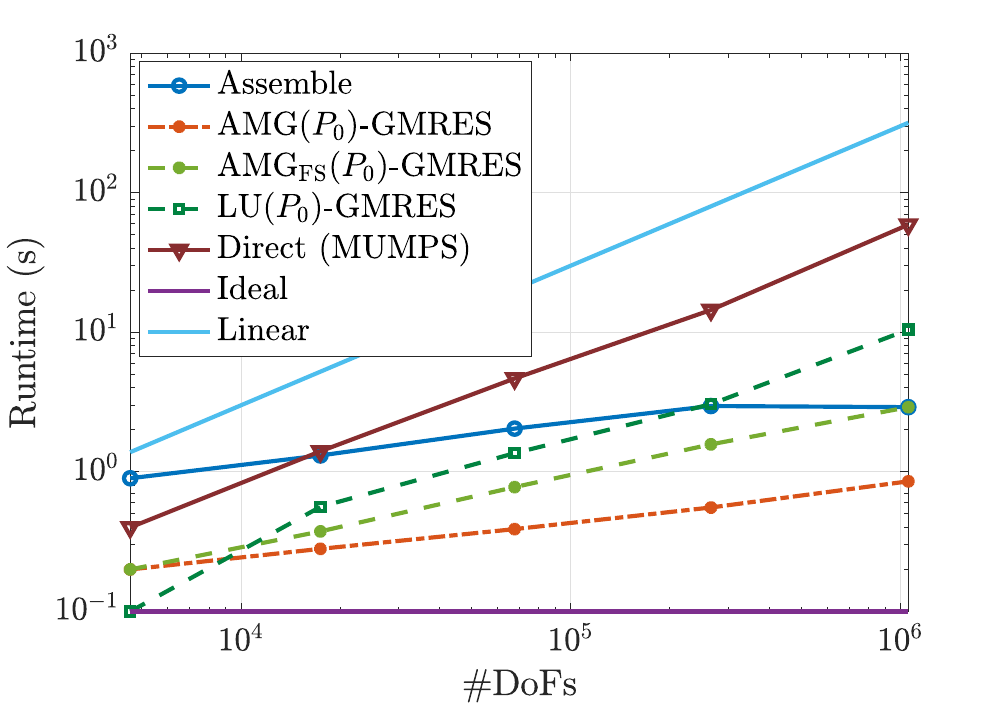}
    \caption{Model A ($d=2$)}
    \label{fig::weak_scaling_square}
    \end{subfigure}    
    \hfill
    \begin{subfigure}[b]{0.479\textwidth}
    \centering
    \includegraphics[width=\textwidth]{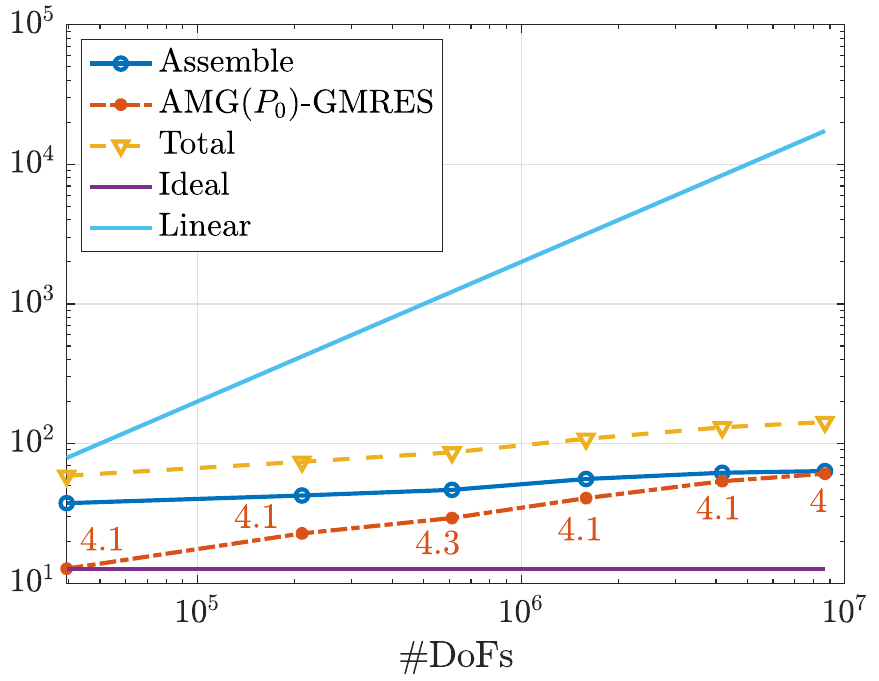}
    \caption{Model A ($d=3$)}
    \label{fig::weak_scaling_cube}
    \end{subfigure}        
    \caption{Parallel performance: weak scaling for Model A ($p = 1$) with (a) $d = 2$, $N_t = 10$, $4.5 \cdot 10^3$ degrees of freedom per core and (b) $d = 3$ and $N_t = 20$, and $4.0 \cdot 10^4$ degrees of freedom per core. Average AMG($P_0$)-GMRES iterations over $n=1,\ldots,N_t$ are displayed as well as the total runtime.}    
    \label{fig::weak_scaling}
\end{figure}

Finally, we investigate the weak scaling (increasing the number of cores and the problem size simultaneously) of the various solution strategies, returning again to Model A (Figure~\ref{fig::weak_scaling}). For small-sized problems (Figure~\ref{fig::weak_scaling}, left), the field-split AMG$_{\text{FS}}(P_0)$-GMRES is the fastest in serial, but, again, the monolithic AMG($P_0$) outperforms all the other approaches in terms of scaling. 
For the moderate-sized problems (Figure~\ref{fig::weak_scaling_cube}, right), near optimal scalability is observed for the assembly, with essentially a flat curve from 64 to 256 cores. As expected, the number of AMG($P_0$)-GMRES iterations stays near constant ($\approx4$), while its runtime increases by an average factor of $1.2\times$ as the number of cores and degrees of freedom doubles. The increase in runtime is connected to communication overhead in matrix-vector products and AMG setup. We refer to \cite{d2021amg} for examples of fine tuning AMG parameters to further improve its weak scalability. 

\section{Discussion, conclusions and outlook}
\label{sec:conclusions}

The availability of dense brain tissue reconstructions at an extreme level of detail poses a unique challenge for mathematical and computational modelling. Here, we presented a scalable finite element discretization and solution algorithm to solve the KNP-EMI equations, describing electrodiffusion, at the level of complexity required for physiologically relevant simulations. In particular, we employ a continuous finite elements discretization of arbitrary order, with possible discontinuities across interfaces, where the KNP-EMI equations are coupled with one or more ODEs, modeling complex membrane dynamics. We describe the algebraic structure of the arising monolithic linear system and present a tailored multilevel solution strategy, since black-box approaches are unfeasible for large problems.

We show robustness of the our proposed solution strategy with respect to discretization parameters, given a near uniform spatial refinement.
Strong and weak scalability of the algorithm are near-optimal, as demonstrated here when running simulations using from 1 to 256 MPI processes for problem sizes ranging from $4.6 \cdot 10^3$ to $10^8$. 
Our motivation for considering this range of parallelism is four-fold: (i) these problem sizes represent the first relevant stage for resolving the brain tissue geometry at the tens of micrometers scale; (ii) any viable solution strategy must successfully address this range; (iii) this is a range of computing resources readily available, in particular 256 represents the default core limit for the Saga computing system; and (iv) above these problem sizes, other aspects of the simulation pipeline such as e.g., mesh generation, labeling, and refinement appear as non-negligible potential bottlenecks. However, our numerical results indicate that the solution algorithm can continue to scale for larger problem sizes and a larger core count. We also remark that while we have presented physiologically relevant scenarios, more complexity could, and should, be considered. For instance, more attention to the connectedness or non-connectedness of the extracellular space and intracellular compartments could yield even higher fidelity representations. Similarly, while we here consider different non-trivial membrane mechanisms for the different cellular types, a higher degree of specificity would be more biologically meaningful, in particular for the neuronal compartments (axons, soma, dendrites). Our simulation code is openly available~\cite{benedusi2024knpemi-zenodo}, and importantly, we envision that the model scenarios presented can be used as stepping stones for future benchmarks and research in this emerging field.  

\section*{Acknowledgments} 
Pietro Benedusi, Ada J. Ellingsrud and Marie E. Rognes acknowledge support from the Research Council of Norway via FRIPRO grant \#324239 (EMIx) and from the national infrastructure for computational science in Norway, Sigma2, via grant \#NN8049K. We thank Dr.~Tom Bartol, Computational Neurobiology Laboratory, Salk Institute for Biological Studies, CA, USA, for providing imaging data and surfaces underlying Model C, as well as inspiring discussions on the topic. We also thank Prof.~Corrado Cal\`i and Prof.~Abdellah Marwan, for sharing surface representations underlying Model B and D and for insightful comments.
Special thanks are also extended to J\o rgen Dokken, Miroslav Kuchta, Francesco Ballarin, Lars Magnus Valnes, Marius Causemann, Rami Masri, Maria Hernandez Mesa, Pasqua D'Ambra, Fabio Durastante, and Salvatore Filippone for valuable comments and pointers.

\section*{Data availability}

This manuscript is accompanied by open software and data sets available at \url{https://doi.org/10.5281/zenodo.10728462}~\cite{benedusi2024knpemi-zenodo}. 

\section*{Conflict of interest}
The authors declare that they have no conflict of interest.

\bibliographystyle{siamplain}
\bibliography{references}

\begin{thebibliography}{10}

\bibitem{abdellah2023ultraliser}
{\sc M.~Abdellah, J.~J.~G. Cantero, N.~R. Guerrero, A.~Foni, J.~S. Coggan,
  C.~Cal{\`\i}, M.~Agus, E.~Zisis, D.~Keller, M.~Hadwiger, et~al.}, {\em
  Ultraliser: a framework for creating multiscale, high-fidelity and
  geometrically realistic {3D} models for in silico neuroscience}, Briefings in
  Bioinformatics, 24 (2023), p.~bbac491.

\bibitem{marwan_abdellah_2022_7105941}
{\sc M.~Abdellah, J.~J.~G. Cantero, N.~R. Guerrero, A.~Foni, J.~S. Coggan,
  C.~Calì, M.~Agus, E.~Zisis, D.~Keller, M.~Hadwiger, P.~J. Magistretti,
  H.~Markram, and F.~Schürmann}, {\em {Ultraliser: a framework for creating
  multiscale, high-fidelity and geometrically realistic 3D models for in silico
  neuroscience}}, Sept. 2022, \url{https://doi.org/10.5281/zenodo.7105941},
  \url{https://doi.org/10.5281/zenodo.7105941}.

\bibitem{agudelo2013computationally}
{\sc A.~Agudelo-Toro and A.~Neef}, {\em Computationally efficient simulation of
  electrical activity at cell membranes interacting with self-generated and
  externally imposed electric fields}, Journal of neural engineering, 10
  (2013), p.~026019.

\bibitem{aitken1986sources}
{\sc P.~Aitken and G.~Somjen}, {\em The sources of extracellular potassium
  accumulation in the {CA1} region of hippocampal slices}, Brain research, 369
  (1986), pp.~163--167.

\bibitem{alnaes2015fenics}
{\sc M.~Aln{\ae}s, J.~Blechta, J.~Hake, A.~Johansson, B.~Kehlet, A.~Logg,
  C.~Richardson, J.~Ring, M.~E. Rognes, and G.~N. Wells}, {\em The {FEniCS}
  project version 1.5}, Archive of Numerical Software, 3 (2015).

\bibitem{armbruster2022neuronal}
{\sc M.~Armbruster, S.~Naskar, J.~P. Garcia, M.~Sommer, E.~Kim, Y.~Adam, P.~G.
  Haydon, E.~S. Boyden, A.~E. Cohen, and C.~G. Dulla}, {\em Neuronal activity
  drives pathway-specific depolarization of peripheral astrocyte processes},
  Nature neuroscience, 25 (2022), pp.~607--616.

\bibitem{petsc-web-page}
{\sc S.~Balay, S.~Abhyankar, M.~F. Adams, S.~Benson, J.~Brown, P.~Brune,
  K.~Buschelman, E.~M. Constantinescu, L.~Dalcin, A.~Dener, V.~Eijkhout,
  J.~Faibussowitsch, W.~D. Gropp, V.~Hapla, T.~Isaac, P.~Jolivet, D.~Karpeev,
  D.~Kaushik, M.~G. Knepley, F.~Kong, S.~Kruger, D.~A. May, L.~C. McInnes,
  R.~T. Mills, L.~Mitchell, T.~Munson, J.~E. Roman, K.~Rupp, P.~Sanan,
  J.~Sarich, B.~F. Smith, S.~Zampini, H.~Zhang, H.~Zhang, and J.~Zhang}, {\em
  {PETS}c {W}eb page}.
\newblock \url{https://petsc.org/}, 2023, \url{https://petsc.org/}.

\bibitem{multiphenics}
{\sc F.~Ballarin}, {\em multiphenics}, 2024,
  \url{https://multiphenics.github.io/}.

\bibitem{bartol2015computational}
{\sc T.~M. Bartol, D.~X. Keller, J.~P. Kinney, C.~L. Bajaj, K.~M. Harris, T.~J.
  Sejnowski, and M.~B. Kennedy}, {\em Computational reconstitution of spine
  calcium transients from individual proteins}, Frontiers in synaptic
  neuroscience, 7 (2015), p.~17.

\bibitem{benedusi2021fast}
{\sc P.~Benedusi, P.~Ferrari, C.~Garoni, R.~Krause, and S.~Serra-Capizzano},
  {\em Fast parallel solver for the space-time {IgA-DG} discretization of the
  diffusion equation}, Journal of Scientific Computing, 89 (2021), p.~20.

\bibitem{benedusi2024EMI}
{\sc P.~Benedusi, P.~Ferrari, M.~Rognes, and S.~Serra-Capizzano}, {\em Modeling
  excitable cells with the {EMI} equations: spectral analysis and iterative
  solution strategy}, Journal of Scientific Computing, 98 (2024), p.~58.

\bibitem{benedusi2024knpemi-zenodo}
{\sc P.~Benedusi and M.~E. Rognes}, {\em Scalable approximation and solvers for
  ionic electrodiffusion in cellular geometries}, 2024,
  \url{https://doi.org/https://doi.org/10.5281/zenodo.10728462}.

\bibitem{berre2023cut}
{\sc N.~Berre, M.~E. Rognes, and A.~Massing}, {\em Cut finite element
  discretizations of cell-by-cell {EMI} electrophysiology models}, arXiv
  preprint arXiv:2306.03001,  (2023).

\bibitem{boron2016medical}
{\sc W.~F. Boron and E.~L. Boulpaep}, {\em Medical physiology}, Elsevier Health
  Sciences, 2016.

\bibitem{buccino2021improving}
{\sc A.~P. Buccino, M.~Kuchta, J.~Schreiner, and K.-A. Mardal}, {\em Improving
  neural simulations with the {EMI} model}, in Modeling Excitable Tissue,
  Springer, 2021, pp.~87--98.

\bibitem{cali20193d}
{\sc C.~Cal{\`\i}, M.~Agus, K.~Kare, D.~J. Boges, H.~Lehv{\"a}slaiho,
  M.~Hadwiger, and P.~J. Magistretti}, {\em {3D} cellular reconstruction of
  cortical glia and parenchymal morphometric analysis from serial block-face
  electron microscopy of juvenile rat}, Progress in neurobiology, 183 (2019),
  p.~101696.

\bibitem{microns2021functional}
{\sc M.~Consortium, J.~A. Bae, M.~Baptiste, C.~A. Bishop, A.~L. Bodor,
  D.~Brittain, J.~Buchanan, D.~J. Bumbarger, M.~A. Castro, B.~Celii, et~al.},
  {\em Functional connectomics spanning multiple areas of mouse visual cortex},
  BioRxiv,  (2021), pp.~2021--07.

\bibitem{d2021amg}
{\sc P.~D'Ambra, F.~Durastante, and S.~Filippone}, {\em {AMG} preconditioners
  for linear solvers towards extreme scale}, SIAM Journal on Scientific
  Computing, 43 (2021), pp.~S679--S703.

\bibitem{dickinson2011electroneutrality}
{\sc E.~J. Dickinson, J.~G. Limon-Petersen, and R.~G. Compton}, {\em The
  electroneutrality approximation in electrochemistry}, Journal of Solid State
  Electrochemistry, 15 (2011), pp.~1335--1345.

\bibitem{dietz2023local}
{\sc A.~G. Dietz, P.~Weikop, N.~Hauglund, M.~Andersen, N.~C. Petersen, L.~Rose,
  H.~Hirase, and M.~Nedergaard}, {\em Local extracellular {K+} in cortex
  regulates norepinephrine levels, network state, and behavioral output},
  Proceedings of the National Academy of Sciences, 120 (2023), p.~e2305071120.

\bibitem{ellingsrud2021accurate}
{\sc A.~J. Ellingsrud, N.~Boull{\'e}, P.~E. Farrell, and M.~E. Rognes}, {\em
  Accurate numerical simulation of electrodiffusion and water movement in brain
  tissue}, Mathematical Medicine and Biology: A Journal of the IMA, 38 (2021),
  pp.~516--551.

\bibitem{ellingsrud2020finite}
{\sc A.~J. Ellingsrud, A.~Solbr{\aa}, G.~T. Einevoll, G.~Halnes, and M.~E.
  Rognes}, {\em Finite element simulation of ionic electrodiffusion in cellular
  geometries}, Frontiers in Neuroinformatics, 14 (2020), p.~11.

\bibitem{falgout2002hypre}
{\sc R.~D. Falgout and U.~M. Yang}, {\em hypre: A library of high performance
  preconditioners}, in Computational Science—ICCS 2002: International
  Conference Amsterdam, The Netherlands, April 21--24, 2002 Proceedings, Part
  III, Springer, 2002, pp.~632--641.

\bibitem{farina2021cut}
{\sc S.~Farina, S.~Claus, J.~S. Hale, A.~Skupin, and S.~P. Bordas}, {\em A cut
  finite element method for spatially resolved energy metabolism models in
  complex neuro-cell morphologies with minimal remeshing}, Advanced Modeling
  and Simulation in Engineering Sciences, 8 (2021), pp.~1--32.

\bibitem{halnes2013electrodiffusive}
{\sc G.~Halnes, I.~{\O}stby, K.~H. Pettersen, S.~W. Omholt, and G.~T.
  Einevoll}, {\em Electrodiffusive model for astrocytic and neuronal ion
  concentration dynamics}, PLoS computational biology, 9 (2013), p.~e1003386.

\bibitem{hille2001ion}
{\sc B.~Hille et~al.}, {\em {Ion channels of excitable membranes}}, vol.~507,
  Sinauer Sunderland, MA, 2001.

\bibitem{hodgkin1952quantitative}
{\sc A.~L. Hodgkin and A.~F. Huxley}, {\em A quantitative description of
  membrane current and its application to conduction and excitation in nerve},
  The Journal of physiology, 117 (1952), p.~500.

\bibitem{ftet}
{\sc Y.~Hu, T.~Schneider, B.~Wang, D.~Zorin, and D.~Panozzo}, {\em Fast
  tetrahedral meshing in the wild}, ACM Trans. Graph., 39 (2020),
  \url{https://doi.org/10.1145/3386569.3392385},
  \url{https://doi.org/10.1145/3386569.3392385}.

\bibitem{huynh2023convergence}
{\sc N.~M.~M. Huynh, F.~Chegini, L.~F. Pavarino, M.~Weiser, and S.~Scacchi},
  {\em Convergence analysis of {BDDC} preconditioners for composite dg
  discretizations of the cardiac cell-by-cell model}, SIAM Journal on
  Scientific Computing, 45 (2023), pp.~A2836--A2857.

\bibitem{jaeger2022arrhythmogenic}
{\sc K.~H. J{\ae}ger, A.~G. Edwards, W.~R. Giles, and A.~Tveito}, {\em
  Arrhythmogenic influence of mutations in a myocyte-based computational model
  of the pulmonary vein sleeve}, Scientific Reports, 12 (2022), p.~7040.

\bibitem{jaeger2021efficient}
{\sc K.~H. J{\ae}ger, K.~G. Hustad, X.~Cai, and A.~Tveito}, {\em Efficient
  numerical solution of the {EMI} model representing the extracellular space
  (e), cell membrane (m) and intracellular space (i) of a collection of cardiac
  cells}, Frontiers in Physics, 8 (2021), p.~579461.

\bibitem{kespe2019three}
{\sc M.~Kespe, S.~Cernak, M.~Glei{\ss}, S.~Hammerich, and H.~Nirschl}, {\em
  Three-dimensional simulation of transport processes within blended electrodes
  on the particle scale}, International Journal of Energy Research, 43 (2019),
  pp.~6762--6778.

\bibitem{kinney2013extracellular}
{\sc J.~P. Kinney, J.~Spacek, T.~M. Bartol, C.~L. Bajaj, K.~M. Harris, and
  T.~J. Sejnowski}, {\em Extracellular sheets and tunnels modulate glutamate
  diffusion in hippocampal neuropil}, Journal of Comparative Neurology, 521
  (2013), pp.~448--464.

\bibitem{kuchta2021solving}
{\sc M.~Kuchta, K.-A. Mardal, and M.~E. Rognes}, {\em Solving the emi equations
  using finite element methods}, Modeling Excitable Tissue: The EMI Framework,
  (2021), pp.~56--69.

\bibitem{lagache2019electrodiffusion}
{\sc T.~Lagache, K.~Jayant, and R.~Yuste}, {\em Electrodiffusion models of
  synaptic potentials in dendritic spines}, Journal of computational
  neuroscience, 47 (2019), pp.~77--89.

\bibitem{lee20203d}
{\sc C.~T. Lee, J.~G. Laughlin, N.~Angliviel~de La~Beaumelle, R.~E. Amaro,
  J.~A. McCammon, R.~Ramamoorthi, M.~Holst, and P.~Rangamani}, {\em {3D} mesh
  processing using {GAMer} 2 to enable reaction-diffusion simulations in
  realistic cellular geometries}, PLoS computational biology, 16 (2020),
  p.~e1007756.

\bibitem{logg2012automated}
{\sc A.~Logg, K.-A. Mardal, and G.~Wells}, {\em Automated solution of
  differential equations by the finite element method: The FEniCS book},
  vol.~84, Springer Science \& Business Media, 2012.

\bibitem{mori2015multidomain}
{\sc Y.~Mori}, {\em A multidomain model for ionic electrodiffusion and osmosis
  with an application to cortical spreading depression}, Physica D: Nonlinear
  Phenomena, 308 (2015), pp.~94--108.

\bibitem{mori2009numerical}
{\sc Y.~Mori and C.~Peskin}, {\em A numerical method for cellular
  electrophysiology based on the electrodiffusion equations with internal
  boundary conditions at membranes}, Communications in Applied Mathematics and
  Computational Science, 4 (2009), pp.~85--134,
  \url{http://dx.doi.org/10.2140/camcos.2009.4.85}.

\bibitem{motta2019dense}
{\sc A.~Motta, M.~Berning, K.~M. Boergens, B.~Staffler, M.~Beining, S.~Loomba,
  P.~Hennig, H.~Wissler, and M.~Helmstaedter}, {\em Dense connectomic
  reconstruction in layer 4 of the somatosensory cortex}, Science, 366 (2019).

\bibitem{nicholson1978calcium}
{\sc C.~Nicholson, G.~Ten~Bruggencate, H.~Stockle, and R.~Steinberg}, {\em
  Calcium and potassium changes in extracellular microenvironment of cat
  cerebellar cortex}, Journal of neurophysiology, 41 (1978), pp.~1026--1039.

\bibitem{niederer2013regulation}
{\sc S.~Niederer}, {\em Regulation of ion gradients across myocardial ischemic
  border zones: a biophysical modelling analysis}, PloS one, 8 (2013),
  p.~e60323.

\bibitem{nordentoft2023local}
{\sc M.~S. Nordentoft, A.~Papoutsi, N.~Takahashi, M.~S. Heltberg, M.~H. Jensen,
  and R.~N. Rasmussen}, {\em Local changes in potassium ions modulate dendritic
  integration}, bioRxiv,  (2023), pp.~2023--05.

\bibitem{pods2013electrodiffusion}
{\sc J.~Pods, J.~Sch{\"o}nke, and P.~Bastian}, {\em {Electrodiffusion models of
  neurons and extracellular space using the Poisson-Nernst-Planck
  equations---numerical simulation of the intra-and extracellular potential for
  an axon model}}, Biophysical journal, 105 (2013), pp.~242--254.

\bibitem{rasmussen2020interstitial}
{\sc R.~Rasmussen, J.~O’Donnell, F.~Ding, and M.~Nedergaard}, {\em
  Interstitial ions: A key regulator of state-dependent neural activity?},
  Progress in Neurobiology,  (2020), p.~101802.

\bibitem{rosilho2024boundary}
{\sc G.~Rosilho~de Souza, R.~Krause, S.~Pezzuto, et~al.}, {\em Boundary
  integral formulation of the cell-by-cell model of cardiac electrophysiology},
  Engineering Analysis with Boundary Elements, 158 (2024), pp.~239--251.

\bibitem{roy2023scalable}
{\sc T.~Roy, J.~Andrej, and V.~A. Beck}, {\em A scalable {DG} solver for the
  electroneutral {Nernst-Planck} equations}, Journal of Computational Physics,
  475 (2023), p.~111859.

\bibitem{saetra2021electrodiffusive}
{\sc M.~J. S{\ae}tra, G.~T. Einevoll, and G.~Halnes}, {\em An electrodiffusive
  neuron-extracellular-glia model for exploring the genesis of slow potentials
  in the brain}, PLoS Computational Biology, 17 (2021), p.~e1008143.

\bibitem{saetra2023neural}
{\sc M.~J. S{\ae}tra, A.~J. Ellingsrud, and M.~E. Rognes}, {\em Neural activity
  induces strongly coupled electro-chemo-mechanical interactions and fluid flow
  in astrocyte networks and extracellular space--a computational study}, PLOS
  Computational Biology,  (2023),
  \url{https://doi.org/10.1371/journal.pcbi.1010996}.

\bibitem{shirinpour2021multi}
{\sc S.~Shirinpour, N.~Hananeia, J.~Rosado, H.~Tran, C.~Galanis, A.~Vlachos,
  P.~Jedlicka, G.~Queisser, and A.~Opitz}, {\em Multi-scale modeling toolbox
  for single neuron and subcellular activity under transcranial magnetic
  stimulation}, Brain stimulation, 14 (2021), pp.~1470--1482.

\bibitem{solbraa2018kirchhoff}
{\sc A.~Solbr{\aa}, A.~W. Bergersen, J.~van~den Brink, A.~Malthe-S{\o}renssen,
  G.~T. Einevoll, and G.~Halnes}, {\em A {Kirchhoff-Nernst-Planck} framework
  for modeling large scale extracellular electrodiffusion surrounding
  morphologically detailed neurons}, PLoS computational biology, 14 (2018),
  p.~e1006510.

\bibitem{sterratt2023principles}
{\sc D.~Sterratt, B.~Graham, A.~Gillies, G.~Einevoll, and D.~Willshaw}, {\em
  Principles of computational modelling in neuroscience}, Cambridge University
  Press, 2023, \url{https://doi.org/10.1017/CBO9780511975899}.

\bibitem{sundnes2006computational}
{\sc J.~Sundnes, B.~F. Nielsen, K.~A. Mardal, X.~Cai, G.~T. Lines, and
  A.~Tveito}, {\em On the computational complexity of the bidomain and the
  monodomain models of electrophysiology}, Annals of biomedical engineering, 34
  (2006), pp.~1088--1097.

\bibitem{tveito2017cell}
{\sc A.~Tveito, K.~H. J{\ae}ger, M.~Kuchta, K.-A. Mardal, and M.~E. Rognes},
  {\em A cell-based framework for numerical modeling of electrical conduction
  in cardiac tissue}, Frontiers in Physics, 5 (2017), p.~48.

\bibitem{tveito2021modeling}
{\sc A.~Tveito, K.-A. Mardal, and M.~E. Rognes}, {\em Modeling excitable
  tissue:~the {EMI} framework}, Springer, 2021.

\bibitem{untiet2023astrocytic}
{\sc V.~Untiet, F.~R. Beinlich, P.~Kusk, N.~Kang, A.~Ladr{\'o}n-de Guevara,
  W.~Song, C.~Kjaerby, M.~Andersen, N.~Hauglund, Z.~Bojarowska, et~al.}, {\em
  Astrocytic chloride is brain state dependent and modulates inhibitory
  neurotransmission in mice}, Nature Communications, 14 (2023), p.~1871.

\bibitem{vossen2007modeling}
{\sc C.~Vo{\ss}en, J.~P. Eberhard, and G.~Wittum}, {\em Modeling and simulation
  for three-dimensional signal propagation in passive dendrites}, Computing and
  Visualization in Science, 10 (2007), pp.~107--121.

\bibitem{xylouris2010three}
{\sc K.~Xylouris, G.~Queisser, and G.~Wittum}, {\em A three-dimensional
  mathematical model of active signal processing in axons}, Computing and
  visualization in science, 13 (2010), pp.~409--418.

\bibitem{xylouris2015three}
{\sc K.~Xylouris and G.~Wittum}, {\em A three-dimensional mathematical model
  for the signal propagation on a neuron's membrane}, Frontiers in
  computational neuroscience, 9 (2015), p.~94.

\bibitem{yang2022carbon}
{\sc H.~Yang, Z.~Qian, J.~Wang, J.~Feng, C.~Tang, L.~Wang, Y.~Guo, Z.~Liu,
  Y.~Yang, K.~Zhang, et~al.}, {\em Carbon nanotube array-based flexible
  multifunctional electrodes to record electrophysiology and ions on the
  cerebral cortex in real time}, Advanced Functional Materials, 32 (2022),
  p.~2204794.

\bibitem{yao2017numerical}
{\sc L.~Yao and Y.~Mori}, {\em A numerical method for osmotic water flow and
  solute diffusion with deformable membrane boundaries in two spatial
  dimension}, Journal of Computational Physics, 350 (2017), pp.~728--746.

\bibitem{zisis2021architecture}
{\sc E.~Zisis, D.~Keller, L.~Kanari, A.~Arnaudon, M.~Gevaert, T.~Delemontex,
  B.~Coste, A.~Foni, A.~Marwan, C.~Cal{\`\i}, et~al.}, {\em Architecture of the
  neuro-glia-vascular system}, bioRxiv 2021.01.19.427241,  (2021).

\end{thebibliography}

\end{document}